\documentclass{amsart}
\usepackage{hyperref}
\usepackage{amssymb}
\usepackage[dvips]{color}
\usepackage{graphicx}
\usepackage{epsfig}
\usepackage[all]{xy}

\newlength{\cellsz}
\newcounter{cellsize}
\newcommand{\setcellsize}[1]{%
  \setcounter{cellsize}{#1}%
  \setlength{\cellsz}{\value{cellsize}\unitlength}}%
\setcellsize{13}%
\newcommand\cellify[1]{\def\thearg{#1}\def\nothing{}%
\ifx\thearg\nothing \vrule width0pt height\cellsz depth0pt\else
\hbox to 0pt{{\begin{picture}(\value{cellsize},\value{cellsize})
  \put(0,0){\line(1,0){\value{cellsize}}}
  \put(0,0){\line(0,1){\value{cellsize}}}
  \put(\value{cellsize},0){\line(0,1){\value{cellsize}}}
  \put(0,\value{cellsize}){\line(1,0){\value{cellsize}}} \end{picture} \hss}}\fi%
\vbox to \cellsz{ \vss \hbox to \cellsz{\hss$#1$\hss} \vss}}
\newcommand\tableau[1]{\vcenter{\vbox{\let\\\cr
\baselineskip -16000pt \lineskiplimit 16000pt \lineskip 0pt
\ialign{&\cellify{##}\cr#1\crcr}}}}
\newcommand\tabl[1]{\vtop{\let\\\cr
\baselineskip -16000pt \lineskiplimit 16000pt \lineskip 0pt
\ialign{&\cellify{##}\cr#1\crcr}}}

\def\sh{\operatorname{sh}}
\def\st{\operatorname{st}}
\def\type{\operatorname{type}}

\def\m{\mathfrak m}
\def\M{\mathfrak M}
\def\MPR{{\rm MR}}
\def\mMPR{{\mathfrak{m}}{\rm MR}}
\def\MMPR{{\mathfrak{M}}{\rm MR}}
\def\MNSym{{\mathfrak{M}}{\rm NSym}}
\def\MSym{{\mathfrak{M}}{\rm Sym}}
\def\bS{S^{\mathfrak{M}}}

\def\A{{\mathcal A}}
\def\Comp{{\rm Comp}}
\def\Z{\mathbb Z}
\def\mQSym{{\mathfrak m}{\rm QSym}}
\def\mSym{{\mathfrak m}{\rm Sym}}
\def\QSym{{\rm QSym}}
\def\NSym{{\rm NSym}}
\def\Sym{{\rm Sym}}

\def\c{{\mathbb C}}

\def\Des{{\rm Des}}

\def\J{\mathcal J}

\def\C{\mathcal C}

\def\P{\mathbb P}
\def\N{\mathbb N}

\def\Z{\mathbb Z}
\def\wt{\mathrm{wt}}

\newtheorem{theorem}{Theorem}[section]
\newtheorem{lemma}[theorem]{Lemma}

\newtheorem{prop}[theorem]{Proposition}
\newtheorem{corollary}[theorem]{Corollary}

\theoremstyle{definition}
\newtheorem{example}[theorem]{Example}
\newtheorem{problem}[theorem]{Problem}
\newtheorem{definition}[theorem]{Definition}

\theoremstyle{remark}

\newtheorem{remark}[theorem]{Remark}

\def\N{\mathbb N}
\def\L{\tilde L}
\def\R{\tilde R}
\def\K{\tilde K}
\def\Z{\mathbb Z}
\def\PP{\tilde {\mathbb P}}
\def\AA{ \tilde{\mathcal A}}
\def\A{\mathfrak A}

\newcommand{\ip}[1]{\langle #1 \rangle}
\def\Gr{{\rm Gr}}
\def\O{{\mathcal O}}
\def\tg{{\tilde{g}}}
\def\I{{\mathcal I}}
\def\SC{{\rm SC}}

\begin{document}

\title{Combinatorial Hopf algebras and $K$-homology of Grassmanians}
\author{Thomas Lam and Pavlo Pylyavskyy}
\thanks{T.L. was partially supported by NSF DMS-0600677,
P.P. was supported by R. Stanley's NSF DMS-0604423.}

\address{T.L.: Department of Mathematics, Harvard, Cambridge, MA, 02138}
\email{tfylam@math.harvard.edu}

\address{P.P.: Department of Mathematics, M.I.T., Cambridge, MA, 02139}
\email{pasha@mit.edu}

\begin{abstract}
Motivated by work of Buch on set-valued tableaux in relation to the
$K$-theory of the Grassmannian, we study six combinatorial Hopf
algebras. These Hopf algebras can be thought of as $K$-theoretic
analogues of the by now classical ``square'' of Hopf algebras
consisting of symmetric functions, quasisymmetric functions,
noncommutative symmetric functions and the Malvenuto-Reutenauer Hopf
algebra of permutations.  In addition, we develop a theory of
set-valued $P$-partitions and study three new families of symmetric
functions which are weight generating functions of reverse plane
partitions, weak set-valued tableaux and valued-set tableaux.
\end{abstract}

\maketitle


\section{Introduction}
The Hopf algebra $\Sym$ of symmetric functions \cite{EC2}, the Hopf
algebra $\QSym$ of quasisymmetric functions \cite{Ges}, the Hopf
algebra $\NSym$ of noncommutative symmetric functions \cite{GKLLRT}
and the Malvenuto-Reutenauer Hopf algebra $\MPR$ of permutations
\cite{MR} can be arranged in the following diagram:
\begin{equation}\label{eq:four}
\xymatrix{
\Sym \ar@{<<-}[r] \ar@{-}[d] & \NSym \ar@{^{(}->}[r] \ar@{-}[d] & \MPR \ar@{-}[d]& \\
\Sym \ar@{^{(}->}[r] & \QSym \ar@{<<-}[r] & \MPR & }
\end{equation}
Here $\Sym$ and $\NSym$ are Hopf subalgebras of $\QSym$ and $\MPR$
respectively, while $\Sym$ and $\QSym$ are Hopf quotients of $\NSym$
and $\MPR$ respectively. The vertical lines denote Hopf duality, so
that $\Sym$ and $\MPR$ are self-dual.

Each of the four Hopf algebras above come with a distinguished basis
(see Section~\ref{sec:hopf}).  In the case of the symmetric
functions $\Sym$ this is the basis $\{s_\lambda\}$ of Schur
functions.  It is well known that (besides many other
manifestations) Schur functions represent the Schubert classes
$[X_\lambda]$ in the cohomology $H^*(\Gr(k,\c^n))$ of the
Grassmannians $\Gr(k,\c^n)$ of $k$-planes in $\c^n$.

In \cite{LS} Lascoux and Sch\"utzenberger introduced the {\it
{Grothendieck polynomials}} as representatives of $K$-theory classes
of structure sheaves of Schubert varieties. Fomin and Kirillov in
\cite{FK} studied these from combinatorial point of view. In
particular, they introduced the {\it {stable Grothendieck
polynomials}} $G_\lambda$, which are symmetric power series obtained
as a limit of Grothendieck polynomials. In \cite{B}, Buch gave a
combinatorial expression for stable Grothendieck polynomials as
generating series of {\it {set-valued tableaux}}.  These symmetric
functions $G_\lambda$ play the role of Schur functions in the
$K$-theory $K^\circ(\Gr(k,n))$ of Grassmannians; roughly speaking
$G_\lambda$ represents the class of the structure sheaf of a
Schubert variety.  Buch studies a bialgebra $\Gamma$ spanned by the
stable Grothendieck polynomials.  Taking the completion of the
bialgebra $\Gamma$, one can define a Hopf algebra which we denote
$\mSym$.

Our investigation began with the observation that Buch's definition
of set-valued tableaux can be extended to a definition of set-valued
$P$-partitions, thus allowing one to define a ``$K$-theoretic''
analogue $\mQSym$ of the Hopf algebra $\QSym$ of quasisymmetric
functions. In fact the the entire diagram (\ref{eq:four}) can be
extended to give the following diagram:
\[
\xymatrix{
\MSym \ar@{<<-}[r] \ar@{-}[d] & \MNSym \ar@{^{(}->}[r] \ar@{-}[d] & \MMPR \ar@{-}[d]& \\
\mSym \ar@{^{(}->}[r] & \mQSym \ar@{<<-}[r] & \mMPR & }
\]

Here $\MSym$ and $\mQSym$ are Hopf quotients of $\MNSym$ and $\mMPR$
respectively, while $\mSym$ and $\MNSym$ are Hopf subalgebras of
$\mQSym$ and $\MMPR$ respectively.  The vertical lines denote Hopf
duality.  Note that $\mSym$ and $\mMPR$ are no longer self-dual as
combinatorial Hopf algebras.

Before describing our results in more detail we make some general
remarks by grouping the six Hopf algebras into two groups: the {\it
{$\m$-world}} consisting of $\mSym$, $\mQSym$ and $\mMPR$, and the
{\it {$\M$-world}} consisting of $\MSym$, $\MNSym$ and $\MMPR$. The
stable Grothendieck polynomials $G_\lambda$ are in some sense
deformations of the Schur functions $s_\lambda$; in particular, the
lowest degree component of $G_\lambda$ is equal to $s_\lambda$.  In
the same spirit, we observe in the $\m$-world that
\begin{enumerate}
\item
the classical basis constitute the {\it lowest} degree components of
the new basis;
\item
the product in the distinguished basis is infinite
(with the exception of $\mSym$); both the product and coproduct consist of
classical terms plus terms of {\it {higher}} degree.
\end{enumerate}
In the $\M$-world we have:
\begin{enumerate}
\item
the classical basis constitute the {\it {highest}} degree components
of the new basis;
\item
the product and coproduct are finite and consist of classical terms
plus terms of {\it {lower}} degree.
\end{enumerate}

Besides the study of the Hopf structure of these six Hopf algebras,
our main results also include: in the context of $\mQSym$, a theory
of set-valued $P$-partitions; and in the context of $\MSym$ and
$\mSym$  the study of three new families of symmetric functions.
These symmetric functions are weight generating functions of {\it
weak set-valued tableaux}, {\it valued-set tableaux} and {\it
reverse plane partitions}.

We now describe the structure and results of this paper in more
detail.  In Section~\ref{sec:hopf} we review some standard results
concerning the four (by now classical) Hopf algebras of
(\ref{eq:four}).

\subsection{$\mMPR$}
As a preliminary step, we introduce in
Section~\ref{sec:multiShuffle} the {\it {multi-shuffle bialgebra}}
which comes with {\it {multi-shuffle product}} and {\it {cuut
coproduct}}, but leave open the question of whether an antipode can
be defined.  In Section~\ref{sec:mMPR} we describe the {\it {small
multi-Malvenuto-Reutenauer Hopf algebra}} $\mMPR$ in terms of a
basis $w \in S^\m_\infty$ of {\it $\m$-permutations}. The fact that
$\mMPR$ has an antipode is delayed till Section~\ref{sec:dual}.

\subsection{$\mQSym$}
In Section~\ref{sec:mQSym} we describe the Hopf algebra $\mQSym$ of
{\it {multi quasisymmetric functions}}.  We first define $\mQSym$ as
a Hopf quotient of $\mMPR$.  Next, for a labeled poset $(P,\theta)$
we define {\it set-valued $P$-partitions} and show that the
generating function $\K_{P,\theta}$ is a quasisymmetric function
which expands as a sum of {\it multi-fundamental quasisymmetric
functions} $\L_\alpha$ over a {\it multi-Jordan-Holder set}
$\J(P,\theta)$.  We show that the Hopf-algebra of formal linear
combinations of the $\L_\alpha$ (as $\alpha$ varies over all
compositions) is isomorphic to $\mQSym$.  In addition we study the
transition matrix between $\{\L_\alpha\}$ and the (classical)
fundamental quasisymmetric functions $\{L_\alpha\}$.

\subsection{$\mSym$}
In Section~\ref{sec:mSym} we recall precisely the relationship
(\cite{B}) between $\mSym$ and the $K$-theory of Grassmannians.  We
briefly discuss the Fomin-Greene method (\cite{FG}) for obtaining
stable Grothendieck polynomials from operators which act on the
space of partitions.

\subsection{$\MMPR$}
In Section~\ref{sec:MMPR} we enter the big $\M$-world by describing
the {\it {big Multi-Malvenuto-Reutenauer Hopf algebra}} $\MMPR$ in
terms of a basis $w \in \bS_\infty$ of {\it $\M$-permutations}. We
show that $\MMPR$ is dual to $\mMPR$ and describe an intriguing
partial order on $\M$-permutations, generalizing the usual weak
order of the symmetric group.

\subsection{$\MNSym$}
In Section~\ref{sec:MNSym} we describe and study the Hopf subalgebra
$\MNSym \subset \MMPR$ of {\it {Multi-noncommutative symmetric
functions}} in the basis $\{\R_\alpha\}$ which are analogues of
(noncommutative) ribbon Schur functions.  We show that $\mQSym$ and
$\MNSym$ are Hopf-dual.

\subsection{$\MSym$}
In Section~\ref{sec:MSym} we describe the Hopf algebra $\MSym$ of
{\it {Multi-symmetric functions}}.  As an abstract Hopf algebra
$\MSym$ is isomorphic to $\Sym$, but $\MSym$ is equipped with a
distinguished basis $\{g_{\lambda}\}$ of {\it dual stable
Grothedieck polynomials} which are weight generating functions of
reverse plane partitions. We show that the $g_{\lambda}$'s are
symmetric and Schur positive, and describe an explicit rule for
decomposing them into basis of Schur functions. We show that $\MSym$
and $\mSym$ are Hopf-dual and that $\{g_\lambda\}$ and
$\{G_\lambda\}$ are dual bases.  We make explicit here the relation
between $\MSym$ and the $K$-homology of Grassmannian: the basis
$\{g_\lambda\}$ represent the classes in $K$-homology of the ideal
sheaves of the boundaries of Schubert varieties.  In
Sections~\ref{sec:vst}-\ref{sec:vast} we introduce, again using the
Fomin-Greene method, {\it {weak set-valued tableaux}} and {\it
{valued-set tableaux}}.  The weight generating functions of these
tableaux describe the images of $G_{\lambda}$ and $g_{\lambda}$
under the involution $\omega: \Sym \to \Sym$ of the symmetric
functions which sends the elementary symmetric functions $e_n$ to
the homogeneous symmetric functions $h_n$.

\bigskip

New analogues of the Loday-Ronco Hopf algebra of planar binary trees
\cite{LR} can also be defined in the spirit of this paper.  These
Hopf algebras will be the subject of another article \cite{ALMP}.

\bigskip
\begin{remark}
The reader mostly interested in the stable Grothendieck polynomials
$G_{\lambda}$ and dual stable Grothendieck polynomials $g_{\lambda}$
can safely restrict his or her attention to Sections
\ref{sec:mQSym}, \ref{sec:mSym} and  \ref{sec:MSym}.
\end{remark}
\bigskip

{\bf Acknowledgements.} We thank Anders Buch for answering a
question concerning the $K$-theory of Grassmannians.  We thank
Jean-Christophe Novelli and Maria Ronco for drawing our attention to
the results of~\cite{KLNPS} and~\cite{PR}. We are grateful to Sergey
Fomin for making useful suggestions concerning the presentation of
the material.  We would also like to thank Mark Shimozono and Mike
Zabrocki for making their manuscript \cite{SZ} available to us.

\section{Four combinatorial Hopf algebras} \label{sec:hopf}
In this section we briefly describe the four combinatorial Hopf
algebras, $\MPR$, $\QSym$, $\NSym$, and $\Sym$ which we intend to
generalize (see (\ref{eq:four})).

We begin with some notation concerning compositions.  A {\it
composition} of $n$ is a sequence $\alpha =
(\alpha_1,\alpha_2,\ldots,\alpha_k)$ of positive integers such that
$\alpha_1 + \alpha_2 + \cdots + \alpha_k = n$. We write $|\alpha| =
n$.  Denote the set of compositions of $n$ by $\Comp(n)$. Associated
to a composition $\alpha = (\alpha_1,\alpha_2,\ldots,\alpha_k)$ of
$n$ is a descent subset $D(\alpha) =
\{\alpha_1,\alpha_1+\alpha_2,\ldots,\alpha_1+\alpha_2+\dots +
\alpha_{k-1}\}$ of $[n-1] = \{1,2,\ldots,n-1\}$.  The map $\alpha
\mapsto D(\alpha)$ is a bijection between compositions of $n$ and
subsets of $[n-1]$.  We will denote the inverse map by $\C:
2^{[n-1]} \to \Comp(n)$ so that $\C(D(\alpha))= \alpha$.

Now if $w \in S_n$ is a permutation we let $\Des(w) = \{i \in [n-1]
\mid w_i > w_{i+1}\}$ denote its descent set, and define $\C(w) =
\C(\Des(w))$.  If $\alpha \in \Comp(n)$, we let $w(\alpha)$ denote
any permutation such that $\C(w(\alpha)) = \alpha$.  Similarly, we
define $w(D)$ for a subset $D$ of $[n-1]$.  Note that when we
compare two descent sets, for example $\Des(w)$ and $D(\alpha)$, we
always compare them as subsets.  Thus $\Des(w)$ and $D(\alpha)$ will
never be equal unless $w \in S_{|\alpha|}$.

A partition $\lambda = (\lambda_1 \geq \lambda_2 \geq \cdots \geq
\lambda_l > 0)$ is a decreasing composition.  To a partition
$\lambda$ one may associate its Young diagram, also denoted
$\lambda$, which (in English notation) is drawn as a set of boxes
top-left justified, with $\lambda_i$ boxes in the $i$-th row.  If
$\mu \subset \lambda$ then one obtains a skew Young diagram
$\lambda/\mu$ by taking the set-theoretic difference of $\lambda$
and $\mu$.  We let $\Lambda$ denote the set of all partitions.

\subsection{The Malvenuto-Reutenauer Hopf algebra of permutations}
Malvenuto and Reutenauer~\cite{MR} have defined a Hopf algebra
structure $\MPR$ on the free $\Z$-module spanned by the set of all
permutations $S_\infty = \cup_n S_n$ (by convention $S_0 =
\{\emptyset\}$ contains the empty permutation). For two permutations
$w \in S_n$ and $v \in S_m$ define the {\it shuffle product} $w
\cdot v$ by
$$w \cdot v = w_1 \ldots w_n \times (v_1+n) \ldots (v_m+n),$$ where
$\times$ denotes the usual shuffle of words: for example, $ab \times
cd = abcd + acbd + acdb + cadb + cdab + cabd$.  Now define the {\it
cut coproduct} $\Delta$ by
$$\Delta(w) = \sum_{[u,v]=w} \st(u) \otimes \st(v),$$ where $[u,v]$
denotes the concatenation of $u$ and $v$ and $\st(.)$ is the
standardization operator that replaces any sequence of distinct
integers to the unique permutation with the same set of inversions.
The unit map $\eta: \Z \to \MPR$ is given by $\eta(1) = \emptyset$
while the counit $\varepsilon: \MPR \to \Z$ extracts the coefficient
of $\emptyset$.  With this data $(\cdot, \Delta, \eta, \emptyset)$,
the space $\MPR$ is a bialgebra.  In addition, it has an antipode
(see \cite{AS}) which endows it with the structure of a Hopf
algebra.

The Hopf algebra $\MPR$ is self dual under the map $w \mapsto
(w^{-1})^*$ where $\{w^*\}$ is the basis of $\MPR^*$ dual to
$\{w\}$.

\subsection{The Hopf algebra of quasisymmetric functions}
A formal power series $f = f(x) \in \Z[[x_1,x_2,\ldots]]$ with
bounded degree is called {\it quasisymmetric} if for any
$a_1,a_2,\ldots,a_k \in \P$ we have
$$
\left[x_{i_1}^{a_1}\cdots x_{i_k}^{a_k} \right]f =
\left[x_{j_1}^{a_1}\cdots x_{j_k}^{a_k} \right]f
$$
whenever $i_1 < \cdots < i_k$ and $j_1 <\cdots < j_k$.  Here
$[x_{i_1}^{a_1}\cdots x_{i_k}^{a_k}]f$ denotes the coefficient of
$x_{i_1}^{a_1}\cdots x_{i_k}^{a_k}$ in $f$. Denote by $\QSym \subset
\Z[[x_1,x_2,\ldots]]$ the ring of quasisymmetric functions.

The ring $\QSym$ has a natural coproduct which can be obtained as
follows.  Let $y_1,y_2,\ldots$ be another set of variables and order
$\{x_i \cup y_j\}$ by $x_1 < x_2 < \cdots < y_1 < y_2 < \cdots$.
With this order each $f \in \QSym$ determines $f(x,y)\in
\Z[[x_1,x_2,\ldots, y_1, y_2, \ldots]]$.  One may rewrite $f(x,y)$
as an element of $\Z[[x_1,x_2,\ldots]] \otimes \Z[[y_1,y_2,\ldots]]$
to obtain a coproduct map $\Delta: \QSym \to \QSym \otimes \QSym$.
With this coproduct, $\QSym$ becomes a graded, connected,
commutative but not cocommutative Hopf algebra. We omit the explicit
formula for the antipode.

There are two distinguished $\Z$-bases of $\QSym$, both labeled by
compositions. Let $\alpha$ be a composition of $n$.  The {\it
monomial quasisymmetric function} $M_\alpha$ is given by
\[
M_\alpha = \sum_{i_1 < \cdots <i_k} x_{i_1}^{\alpha_k} \cdots
x_{i_k}^{\alpha_k}.
\]
The {\it fundamental quasisymmetric function} $L_\alpha$ is given by
\[
L_\alpha = \sum_{D(\beta) \supset D(\alpha)} M_\beta =
\sum_{\substack{i_1 \leq i_2 \leq \cdots \leq i_n \\ k \in D(\alpha)
\Rightarrow i_{k+1} > i_k}} x_{i_1} x_{i_2} \cdots x_{i_n},
\]
where the first summation is over compositions $\beta$ satisfying
$|\beta| = |\alpha|$.

The surjective map $w \mapsto L_{\C(w)}$ is a Hopf morphism
(see~\cite{MR}) which exhibits $\QSym$ as a quotient of $\MPR$.

\subsection{The Hopf algebra of noncommutative symmetric functions}
\label{sec:NSym} We refer to \cite{T} for more details concerning
the material of this section.  Let $\NSym$ denote the subspace of
$\MPR$ spanned by the elements
$$
R_{\alpha} = \sum_{\Des(w^{-1})=D(\alpha)} w
$$
for each composition $\alpha$.  It turns out that $\NSym$ is a
graded, connected, cocommutative but not commutative Hopf subalgebra
of $\MPR$. In the basis $\{R_\alpha\}$, the multiplication can be
written as
$$
R_\alpha \,R_\beta = R_{\alpha \vartriangleright \beta} + R_{\alpha \vartriangleleft \beta},
$$
where for $\alpha = (\alpha_1,\ldots,\alpha_k)$ and $\beta =
(\beta_1,\ldots,\beta_l)$ we have $\alpha \vartriangleright \beta =
(\alpha_1,\ldots,\alpha_k+\beta_1,\ldots,\beta_l)$ and $\alpha \vartriangleleft
\beta = (\alpha_1,\ldots,\alpha_k,\beta_1,\ldots,\beta_l)$.

\begin{remark}
Note that in \cite{T} the notation $\alpha \cdot \beta$ is used for
what we call here  $\alpha \vartriangleleft \beta$. We however
change the notation so that $\alpha \cdot \beta$ is saved for a
different operation to be defined later.
\end{remark}

In fact $\NSym$ is a free (noncommutative) algebra generated by
symbols $\{S_i \mid i \geq 1\}$.  Given a composition $\alpha$ we
let $S_\alpha = S_{\alpha_1}S_{\alpha_2} \cdots S_{\alpha_k}$.  Then
$$
S_\alpha = \sum_{D(\beta) \subset D(\alpha)} R_\beta.
$$
In terms of the generators $\{S_i \mid i \geq 1\}$, we have $\Delta
S_i = \sum_{k = 0}^{i} S_k \otimes S_{i-k}$.

The two Hopf algebras $\NSym$ and $\QSym$ are Hopf dual with
$\{R_{\alpha}\}$ and $\{L_{\alpha}\}$ forming dual bases.

\subsection{The Hopf algebra of (commutative) symmetric functions} \label{sec:comsym}
A formal power series $f = f(x) \in \Z[[x_1,x_2,\ldots]]$ with
bounded degree is called {\it symmetric} if for any
$a_1,a_2,\ldots,a_k \in \P$ we have
\[
\left[x_{i_1}^{a_1}\cdots x_{i_k}^{a_k} \right]f =
\left[x_{j_1}^{a_1}\cdots x_{j_k}^{a_k} \right]f
\]
whenever $i_1, \ldots, i_k$ are all distinct and $j_1, \ldots,j_k$
are all distinct.  Denote by $\Sym \subset \Z[[x_1,x_2,\ldots]]$ the
algebra of symmetric functions.  Every symmetric function is
quasisymmetric and in fact $\Sym$ is a commutative and cocommutative
Hopf subalgebra of $\QSym$.

The Hopf algebra of symmetric functions has a distinguished basis
$\{s_\lambda \mid \lambda \in \Lambda\}$ of Schur functions indexed
by the set of all partitions.  The Schur function $s_\lambda$ is the
weight generating function of semistandard tableaux with shape
$\lambda$.  The {\it Hall inner product} of $\Sym$ is defined by
$\ip{s_\lambda,s_\mu} = \delta_{\lambda\mu}$.  With this inner
product, $\Sym$ is a self-dual Hopf algebra.  The product and
coproduct structure constants of $\Sym$ are both given by the {\it
Littlewood-Richardson coefficients} $c^\nu_{\lambda\mu} \in \Z$:
$$
s_\lambda s_\mu = \sum_\nu c^\nu_{\lambda\mu} s_nu \ \ \ \ \
\Delta(s_\lambda) = \sum_{\mu,\nu} c^\lambda_{\mu\nu} s_\mu \otimes
s_\nu.
$$

A skew shape $\lambda/\mu$ is a ribbon if it is connected and
contains no $2 \times 2$ square.  Given a composition $\alpha =
(\alpha_1,\alpha_2,\ldots,\alpha_k)$, there is a ribbon $r_\alpha =
\lambda/\mu$ with $\alpha_k$ boxes in its first row, $\alpha_{k-1}$
boxes in its second row and so on.  The map $R_\alpha \mapsto
s_{r_\alpha}$ expresses $\Sym$ as the commutative quotient of
$\NSym$.


\section{Multi-shuffle algebra}
\label{sec:multiShuffle}

In the Sections~\ref{sec:mMPR}-\ref{sec:mSym} we will define and
study three Hopf algebras: the small multi-Malvenuto-Reutenauer
algebra $\mMPR$, the algebra of multi-quasisymmetric functions
$\mQSym$ and the algebra of multi-symmetric functions $\mSym$.  The
Hopf algebra $\mQSym$ is a Hopf quotient of $\mMPR$ while $\mSym$ is
a Hopf subalgebra of $\mQSym$.  Now, as a preliminary step, we
define the multi-shuffle algebra.

In the following we will be dealing with free $\Z$-modules $M$ which
are the sets of {\it arbitrary} $\Z$-linear combinations of a
countable set $S$.  If $N$ is another $\Z$-module, we will say that
a linear function $M \to N$ is {\it continuous} if it respects
arbitrary linear combinations of elements of $S$ (not just finite
linear combinations).

Let $\A$ denote an alphabet, $\A^*$ denote
the set of (possibly empty) words with letters from $\A$ and let $\m
S[\A^*] = \prod_{a \in \A^*} \Z a$ denote the $\Z$-module of
(infinite) $\Z$-linear combinations words from $\A^*$.

 Let $u = u_1 u_2 \cdots u_k \in \A^*$ be a word.  Call  $w = w_1 w_2 \cdots w_{m}$ a {\it
{multiword}} of $u$ if there is a surjective and non-decreasing map
$t: [m] \longrightarrow [k]$, so that $w_j = u_{t(j)}$.  Let $u =
u_1 u_2 \cdots u_k$ and $v = v_1 v_2 \cdots v_l$ be two words, and
assume that all letters $v_i$ and $u_j$ are distinct.  Then a word
$w = w_1 w_2 \cdots w_{m}$ is a {\it {multishuffle}} of $u$ and $v$
if

\begin{enumerate}
 \item neighboring letters of $w$ are distinct, that is for any $i$ we have $w_i \not = w_{i+1}$;
 \item when restricted to alphabets $\{v_i\}$ and $\{u_j\}$ word $w$ becomes a multiword of $v$ and $u$ correspondingly.
\end{enumerate}

We denote by $u \star v$ the sum in $\m S[\A^*]$ of all
multishuffles of $u$ and $v$.  Now suppose that $x = x_1 x_2 \cdots
x_k$ and $y = y_1 y_2 \cdots y_l$ are two other words, with possibly
repeated letters.  Then $x \star y$ is obtained from $u \star v$ by
changing each $u_i$ to $x_i$ and each $v_j$ to $y_j$.  It is
possible to attain multiplicities in this way.  For example, we have
$$ab \star a = aba + 2 aab + 2 aaab + 2 aaba + abab + \ldots.$$  We
extend $\star$ to $\m S[\A^*]$ by linearity and continuity (it is a
quick check to verify that this extension is well-defined).
 One can give the following recursive definition of
multishuffle product.

\begin{prop}
Let $v = v_1 v'$ and $u = u_1 u'$ where $v_1,u_1$ are letters and
$v',u'$ are words.  Then
$$v \star u = (v_1 + u_1 v_1 + v_1 u_1 v_1 + \ldots)(v' \star u) +
(u_1 + v_1 u_1 + u_1 v_1 u_1 + \ldots)(v \star u').$$
\end{prop}

We shall also use the same notation for the multiplication map
$\star: \m S[\A^*] \otimes \m S[\A^*] \to \m S[\A^*]$.
Multishuffling is commutative and associative:

\begin{lemma}
For any three words $u,v,x \in \A^*$, we have $u \star v = v \star
u$ and $(u \star v) \star x = u \star (v \star x)$.
\end{lemma}
\begin{proof}
The first statement is immediate from the definition.  For the
second statement, first assume that all the letters in $u,v,x$ are
distinct.  Then both $(u \star v) \star x$ and $u \star (v \star x)$
are equal to the set of words $w$ satisfying (a) for any $i$, we
have $w_i \neq w_{i+1}$, and (b) when restricted to alphabets
$\{v_i\}$, $\{u_j\}$ and $\{x_k\}$, the word $w$ becomes a multiword
of $v$, $u$ and $x$ correspondingly.  The general case also follows
immediately.
\end{proof}

We now define the {\it cuut} coproduct structure on $\m S[\A^*]$, in
analogy with the cut coproduct.  Let $a = a_1 a_2 a_3 \ldots a_n \in
\A^*$. Then define
$$\blacktriangle(a) = \emptyset \otimes a_1 a_2 a_3 \ldots a_n + a_1 \otimes
a_1 a_2 a_3 \ldots a_n + a_1 \otimes a_2 a_3 \ldots a_n +a_1 a_2
\otimes a_2 a_3 \ldots a_n +$$
$$a_1 a_2 \otimes a_3 \ldots a_n + \cdots + a_1 a_2 a_3 \ldots a_n
\otimes a_n + a_1 a_2 a_3 \ldots a_n \otimes \emptyset.$$ When a
letter $a_i$ occurs twice in a term in the above expression we say
that $a_i$ has been ``cut in the middle'' to obtain such a term. For
example, we have $$\blacktriangle(cut) = \emptyset \otimes cut + c
\otimes cut + c \otimes ut + cu \otimes ut + cu \otimes t + cut
\otimes t + cut \otimes \emptyset.$$

We extend $\blacktriangle$ to $\m S[\A^*]$ by linearity and
continuity.  Define the {\it unit} map $\eta: \Z \to \m S[\A^*]$ by
$\eta(n) = n.\emptyset$ and the counit map $\varepsilon: \m S[\A^*]
\to \Z$ by letting $\varepsilon$ take the coefficient of
$\emptyset$.

\begin{theorem}\label{thm:mshuffle}
The space $\m S[\A^*]$ forms a bialgebra with multi-shuffle product
$\star$, cuut coproduct $\blacktriangle$, unit $\eta$ and counit
$\varepsilon$.
\end{theorem}

We call $\m S[\A^*]$ the {\it multi-shuffle algebra}.

\begin{proof} It is easy to verify that $\m S[\mathfrak A^*]$ is
both a unital associative algebra and a counital coassociative
coalgebra.  We now verify the compatibility of $\star$ and
$\blacktriangle$. Let $w$ and $u$ be two words which we assume for
simplicity to have distinct letters. Then $\blacktriangle(w \star
u)$ is a linear combination of all terms $x \otimes y$ such that
either (1) $x$ is a term in $w' \star u'$ and $y$ is a term in $w''
\star u''$, where $w = w'w''$ and $u = u'u''$, or (2) $x$ is a term
in $w'a \star u'$ and $y$ is a term in $aw'' \star u''$ where $w =
w'aw''$ and $u = u'u''$, or (3) $x$ is a term in $w' \star u'b$ and
$y$ is a term in $w'' \star bu''$ where $w = w'w''$ and $u =
u'bu''$, or (4) $x$ is a term in $w'a \star u'b$ and $y$ is a term
in $aw'' \star bu''$ where $w = w'aw''$ and $u = u'bu''$.  Here
$a,b$ are letters while $w',w'',u',u''$ are words. For example, case
(2) or (4) occurs if some letter $a$ in $w$ lies to both sides of
the cutting point or is cut in the middle when we apply
$\blacktriangle$; otherwise case (1) or (3) occurs.  We check that
the same four kind of terms occur in $\blacktriangle(w) \star
\blacktriangle(u)$.
\end{proof}

\begin{problem}
Does $\m S[\A^*]$ have an antipode?
\end{problem}

\section{The small multi-Malvenuto-Reutenauer Hopf algebra}
\label{sec:mMPR}
\begin{definition}
A {\it {small multi-permutation}}, or $\m$-permutation of $[n]$ is a
word $w$ in the alphabet $1, \ldots, n$ such that no two consecutive
letters in $w$ are equal.  The {\it length} $\ell(w)$ of $w$ is the
the number of letters in $w$.
\end{definition}
We denote the set of multi-permutations of $[n]$ by $S^\m_n$ and the
set of all multi-permutations by $S^\m_\infty = \cup_{n \geq 0}
S^\m_n$. By convention $S^\m_0$ contains a single element -- the
empty multi-permutation $\emptyset$. We let $\mMPR = \prod_{w \in
S^\m_\infty} \Z.w$ be the free $\Z$-module of arbitrary $\Z$-linear
combinations of multi-permutations.

Let $w = w_1 \ldots w_k$ and $u = u_1 \ldots u_l$ be two
multi-permutations, and assume that $w \in S^\m_n$, $u \in S^\m_m$.
Define the product $w * u$ of $w$ and $u$ as follows:
$$w * u = w \star (u+n) = w_1 \ldots w_k \star (u_1 +n)
\ldots (u_l+n)$$ where to use the shuffle product $\star$ we treat
$w$ and $u+n$ as words in the alphabet $\N$.  We extend the formula
by linearity and continuity to give a multiplication $*: \mMPR
\otimes \mMPR \to \mMPR$.

Define the standardization operator $\st: \N^* \to S^\m_\infty$ by
sending a word $w$ to the unique $u \in S^\m_\infty$ of the same
length (if it exists) such that $w_i \leq w_j$ if and only if $u_i
\leq u_j$ for each $1 \leq i,j \leq \ell(w)$.  (Recall that $\N^*$
denotes the set of words in the alphabet $\{1,2,3,\ldots\}$.)  We
define the coproduct $\Delta(w)$ by extending via linearity and
continuity the cuut coproduct as follows:
$$\Delta(w) = \st(\blacktriangle w)$$
where we have extended by linearity and continuity the definition of
$\st$.

Define the unit map $\eta: \Z \to \mMPR$ by $\eta(n) = n \cdot
\emptyset$ and the counit map $\varepsilon: \mMPR \to \Z$ by letting
$\varepsilon$ take the coefficient of $\emptyset$.

\begin{theorem}\label{thm:mMPR}
The space $\mMPR$ is a bialgebra with product $*$, coproduct
$\Delta$, unit $\eta$, and counit $\varepsilon$.
\end{theorem}
\begin{proof}
That $*$ is associative and $\Delta$ is coassociative follows from
the corresponding properties in the multi-shuffle algebra.  The
proof that $*$ and $\Delta$ are compatible is essentially the same
as for Theorem~\ref{thm:mshuffle}.  The only observation needed is
that the standardization operator $\st$ can be applied at the very
end of the calculation of $\Delta(w) * \Delta(u)$ instead of
immediately after calculating $\blacktriangle(w)$ and
$\blacktriangle(u)$.
\end{proof}

We shall call the bialgebra $\mMPR$ of Theorem~\ref{thm:mMPR} the
{\it small multi-Malvenuto-Reutenauer bialgebra}.  We shall show
later that $\mMPR$ has an antipode, making it a Hopf algebra.

\begin{remark}
Call an element $w \in S^\m_\infty$ {\it {irreducible}} if it cannot
be written in the form $w = {v / u} = v_1 \ldots v_k (u_1 +n) \ldots
(u_l+n)$ for two smaller $\m$-permutations $v \in S^\m_n$ and $u$.
 The following simple observations say that combinatorially $\mMPR$
is ``free'' over the set of irreducible elements (this statement is
difficult to make precise because $\mMPR$ is a completion):
\begin{enumerate}
\item Every
$\m$-permutation $w \in S^\m_\infty$ can be uniquely written as $w^1
/ w^2 / \cdots / w^k$ where the $w^i \in S^\m_\infty$ are
irreducible.  We say that $w$ is $k$-reducible in this case.
\item
If $w^1,w^2,\ldots,w^k$ are irreducible, the only term in $w^1
* w^2* \cdots * w^k$ which is $k$-reducible is $w^1
/ w^2 / \cdots / w^k$.
\end{enumerate}
\end{remark}

%

\section{Set-valued $P$-partitions and multi-quasisymmetric functions} \label{sec:mQSym}

\subsection{The Hopf-algebra $\mQSym$}
For $w \in S^\m_\infty$ we define the descent set $\Des(w) \subset
[1,\ell(w)-1]$ by
$$
\Des(w) = \{i \in [1,\ell(w)-1] \mid w_i > w_{i+1}\}.
$$
Thus for example $\Des(15132342) = \{2,4,7\}$.  Note that as in
Section~\ref{sec:hopf} by convention descent sets are always
considered as subsets (of $[1,\ell(w)-1]$) so the descent sets of
$w, u \in S^\m_\infty$ can only coincide if $\ell(w) = \ell(u)$. Let
$I \subset \mMPR$ denote the free $\Z$-submodule spanned by the
elements $w - u$ for pairs $w,u$ satisfying $\Des(w) = \Des(u)$.

\begin{lemma}\label{lem:mQSym}
The subspace $I$ is a biideal of $\mMPR$.  In other words, we have
$I * \mMPR \subset I$, $\mMPR * I \subset I$, and $\Delta(I) \subset
I \otimes \mMPR + \mMPR \otimes I$.
\end{lemma}
\begin{proof}
A term $v \in S^\m_\infty$ occurring in the product $w * u$ of $w
\in S^\m_n$ and $u \in S^\m_m$ is determined by the following
information: (1) the location in $v$ of letters in $[1,n]$ (and
hence also $[n+1,n+m]$), (2) the restriction of $v$ to the alphabet
$[1,n]$ (which is a multiword of $w$), and (3) the restriction of
$v$ to the alphabet $[n+1,n+m]$ (which is a multiword of $u$).  If
$w,w' \in S^\m_\infty$ have the same length then there is a
canonical bijection between the multiwords of $w$ and of $w'$.  Thus
if $w,w',u,u' \in S^\m_\infty$ are such that $\ell(w) = \ell(w')$
and $\ell(u) = \ell(u')$ then there is a canonical bijection $\Phi$
between the terms in $w *u$ and those in $w' * u'$.  If in addition,
$w$ and $w'$ have the same descent set and $u$ and $u'$ have the
same descent set then $\Phi$ preserves descent sets.  This proves
that $I$ is an (algebra) ideal.

Now suppose $w, u \in S^\m_\infty$ satisfy $\Des(w) = \Des(u)$ and
thus in particular $\ell(w) = \ell(u)$.  There is a canonical
bijection between the terms of $\Delta(w)$ and those of $\Delta(u)$
so that if $w'\otimes w''$ corresponds to $u' \otimes u''$ then
$\Des(w') = \Des(u')$ and $\Des(w'') = \Des(u'')$. We show that $w'
\otimes w'' - u' \otimes u'' \in I \otimes \mMPR + \mMPR \otimes I$.
This follows from
$$
w' \otimes w'' - u' \otimes u'' = w' \otimes(w'' - u'') + (w'-u')
\otimes u''.
$$
\end{proof}
By Lemma~\ref{lem:mQSym}, the quotient space $\mMPR/I$ is a quotient
Hopf algebra, which we call $\mQSym$.  We will give an explicit
model for $\mQSym$ as an algebra of quasisymmetric functions.  In
\cite{Haz}, Hazewinkel shows that $\QSym$ is a free algebra over
$\Z$.  It would be interesting to investigate this for $\mQSym$ as
well.


\subsection{Posets and $P$-partitions}
We recall the basic definitions concerning
$P$-partitions~\cite{EC2}.  Let $P$ be a finite poset with $n$
elements and $\theta: P \to [n]$ be a bijective labeling of $P$.

\begin{definition}
A {\it $(P, \theta)$-partition} is a map $\sigma: P \to \mathbb P$
such that for each covering relation $s \lessdot t$ in $P$ we have
\begin{align*}
\sigma(s) &\leq \sigma(t) & \mbox{ if $\theta(s) < \theta (t)$,} \\
\sigma(s) &< \sigma(t) & \mbox{ if $\theta(t) < \theta (s)$.}
\end{align*}
\end{definition}

Denote by ${\mathcal{A}}(P,\theta)$ the set of all
$(P,\theta)$-partitions.  If $P$ is finite then one can define the
formal power series $ K_{P,\theta}(x_1,x_2,\ldots) \in
\Z[[x_1,x_2,\ldots]]$ by
\[
K_{P,\theta}(x_1,x_2,\ldots) = \sum_{\sigma \in
{\mathcal{A}}(P,\theta)} x_1^{\# \sigma^{-1}(1)} x_2^{\#
\sigma^{-1}(2)} \cdots.
\]
The composition $\wt(\sigma) = (\# \sigma^{-1}(1), \#
\sigma^{-1}(2), \ldots)$ is called the {\it weight} of $\sigma$.

Recall that a linear extension of $P$ is a bijection $e: P \to
\{1,2,\ldots,n\}$ satisfying $e(x) < e(y)$ if $x < y$ in $P$.  The
Jordan-Holder set $\J(P,\theta)$ of $(P,\theta)$ is the set
$$\{\theta(e^{-1}(1)) \theta(e^{-1}(1)) \cdots \theta(e^{-1}(n))\mid
\text{$e$ is a linear extension of $P$}\} \subset S_n.$$

\begin{theorem}[\cite{Sta}]
\label{thm:Ppart} The generating function $K_{P,\theta}$ is
quasisymmetric. We have $ K_{P,\theta} = \sum_{w \in \J(P,\theta)}
L_{D(w)}$.
\end{theorem}

\subsection{Set-valued $P$-partitions}
In this section we define set-valued $P$-partitions.  Let $\PP$ be
the set of all non-empty finite subsets of $\mathbb P$. For $a \in
\PP$ we define $\min(a)$ and $\max(a)$ to be the minimal and maximal
elements of $a$.  Suppose that $a,b \in \PP$.  Then we say $a$ is
less than $b$ and write $a \leq b$ if and only if $\max(a) \leq
\min(b)$; similarly we have $a < b$ if and only if $\max(a) <
\min(b)$. Note that $a \leq a$ if and only if $|a|=1$, that is $a$
consists of just one element. Let $(P, \theta)$ as before be a poset
with a bijective labeling.

\begin{definition}
A {\it $(P, \theta)$-set-valued partition} is a map $\sigma: P \to
\PP$ such that for each covering relation $s \lessdot t$ in $P$ we
have
\begin{align*}
\sigma(s) &\leq \sigma(t) & \mbox{ if $\theta(s) < \theta (t)$,} \\
\sigma(s) &< \sigma(t) & \mbox{ if $\theta(t) < \theta (s)$.}
\end{align*}
\end{definition}

Denote by $\AA(P,\theta)$ the set of all $(P,\theta)$-set-valued
partitions.  For $\sigma \in \AA(P,\theta)$ denote by $|\sigma| =
\sum_{a \in P} |\sigma(a)|$ the total number of letters used in the
set-valued partition.  For each $i \in \P$ we define $\sigma^{-1}(i)
= \{ x \in P \mid i \in \sigma(x)\}$.  As before, the composition
$\wt(\sigma) = (\# \sigma^{-1}(1), \# \sigma^{-1}(2), \ldots)$ is
called the {\it weight} of $\sigma$. Define the formal power series
$\K_{P,\theta}(x_1,x_2,\ldots) \in \Z[[x_1,x_2,\ldots]]$ by
\[
\K_{P,\theta}(x_1,x_2,\ldots) = \sum_{\sigma \in \AA(P,\theta)}
x_1^{\# \sigma^{-1}(1)} x_2^{\# \sigma^{-1}(2)} \cdots.
\]

It is easy to see that $\K_{P,\theta}$ is always a quasisymmetric
function. In the following example which will be of major importance
for us $\K_{P,\theta}$ happens to be a symmetric function.

\begin{example} \label{ex:s}

Let $P = \lambda$ be the poset of squares in the Young diagram of a
partition $\lambda = (\lambda_1,\ldots,\lambda_l)$. Let $\theta_{s}$
be the labeling of $\lambda$ obtained from the bottom to top
row-reading order; in other words the bottom row of $\lambda$ is
labeled $1,2,\ldots,\lambda_l$, the next row is labeled $\lambda_l +
1,\ldots, \lambda_{l-1} + \lambda_l$ and so on.  Then
$K_{\lambda,\theta_{s}}$ is equal to the Schur function $s_\lambda$,
while $\K_{\lambda,\theta_{s}}$ is (nearly) equal to the {\it
{stable Grothendieck polynomial}} $G_{\lambda}$ studied in~\cite{B}.
We will return to this example in Section~\ref{sec:mSym}.
\end{example}

Let $\alpha \vdash n$ be a composition and let $C$ be a chain $c_1 <
c_2 < \ldots < c_n$ with $n$ elements and $w = w_1w_2\dots w_n \in
S_n$ a permutation of $\{1,2,\ldots,n\}$ such that $\C(w)= \alpha$.
Then $(C,w)$ can be considered a labeled poset, where $w(c_i) =
w_i$. Now define the {\it multi-fundamental quasisymmetric
functions} by $\tilde L_{\alpha} = \K_{C,w}$ (clearly $\K_{C,w}$
depends only on $\alpha$).

A {\it {linear multi-extension}} of $P$ by $[N]$ is a map $e: P \to
2^{\{1,2,\ldots,N\}}$ for some $N \geq n = |P|$ satisfying
\begin{enumerate}
 \item $e(x) < e(y)$ if $x < y$ in $P$,
 \item each $i \in [N]$ is in $e(x)$ for exactly one $x \in P$, and
 \item none of the sets $e(x)$ contains both $i$ and $i+1$ for any $i$.
\end{enumerate}
The {\it {multi-Jordan-Holder set}} $\tilde \J(P,\theta) = \cup_N
\tilde \J_N(P,\theta)$ of $(P,\theta)$ is the union of sets
$$
\tilde \J_N(P,\theta) = \{\theta(e^{-1}(1)) \theta(e^{-1}(2)) \cdots
\theta(e^{-1}(N))\}
$$ where in the above formula $e$ varies over the
set of linear multi-extensions of $P$ by $[N]$.  The set $\tilde
\J_N(P,\theta)$ is a subset of the set $\{w \in S^\m_n \mid \ell(w)
= N\}$ of $\m$-permutations of length $N$ on $n$ letters.  In fact
if $P$ is the antichain with $n$-elements and $\theta$ is any
labeling then we have $\J_N(P,\theta) = \{w \in S^\m_n \mid \ell(w)
= N\}$.  The following result is the set-valued analogue of
Theorem~\ref{thm:Ppart}.

\begin{theorem}
\label{thm:multiPpart} We have $\K_{P,\theta} = \sum_{N \geq n}
\sum_{w \in \tilde \J_N(P,\theta)} {\tilde L_{\C(w)}}$.
\end{theorem}

\begin{proof}
We give an explicit weight-preserving bijection between
$\AA(P,\theta)$ and the set of pairs $(w,\sigma')$ where $ w\in
\tilde \J_N(P,\theta)$ and $\sigma' \in \AA(C,w)$ where $C = (c_1 <
c_2 < \cdots < c_l)$ is a chain with $l = \ell(w)$ elements.  Let
$\sigma \in \AA(P,\theta)$. For each $i$, identify $\sigma^{-1}(i)$
with a subset of $[n]$ via $\theta$ and let $w^{(i)}_\sigma$ denote
the word of length $|\sigma^{-1}(i)|$ obtained by writing the
elements of $\sigma^{-1}(i)$ in increasing order.  Let $w$ be the
unique $\m$-permutation such that $w_\sigma:=w_\sigma^{(1)}
w_\sigma^{(2)} \cdots$ is a multiword of $w$, and we let $t:
\ell(w_\sigma) \to \ell(w)$ denote the associated function as in
Section~\ref{sec:multiShuffle}.  Note that $\sigma^{(r)} =
\emptyset$ for sufficiently large $r$, so that $w_\sigma$ is a
finite word in the alphabet $[n]$ (using all the letters of $[n]$).
Now define $\sigma' \in \AA(C,w)$ by
$$
\sigma'(c_i) = \{r \in \P \mid \text{$w_\sigma^{(r)}$ contributes
letters to $w_\sigma|_{t^{-1}(i)}$}\}
$$
where $w_\sigma|_{t^{-1}(i)}$ is the set of letters in $w_\sigma$ at
the positions in the interval $t^{-1}(i)$.  We claim that this
defines a map $\alpha:\sigma \mapsto (w,\sigma')$ with the required
properties.

First, $w$ is the multiword associated to the linear multi-extension
$e_w$ of $P$ by $\ell(w)$ defined by the condition that $e_w(x)$
contains $j$ if and only if $w_j = \theta(x)$.  It follows from the
definition that this $e_w: P \to 2^{[1,\ell(w)]}$ is a linear
multi-extension.  To check that $\sigma'$ is a set-valued $(C,w)$
partition, we note that $\sigma'(c_i)\leq \sigma'(c_{i+1})$ by
definition, since the function $t$ is non-decreasing.  Furthermore,
if $w_i > w_{i+1}$ then $\sigma'(c_i) < \sigma'(c_{i+1})$ because
each $w_\sigma^{(r)}$ was defined to be increasing.

Finally, the inverse map $\beta: (w,\sigma') \mapsto \sigma$ can be
defined by the formula
$$
\sigma(x) = \bigcup_{j \in e_w(x)} \sigma'(c_j).
$$
That the above union is disjoint follows from the fact that there is
always a descent somewhere between two occurrences of the same
letter in $w$.  That $\sigma$ as defined respects $\theta$ is due to
the fact that $e_w$ is a linear multi-extension.  The equation
$\beta \circ \alpha = {\rm id}$ follows immediately.  For $\alpha
\circ \beta = {\rm id}$, consider a subset $\sigma'(c_j) \subset
\sigma(x)$.  One checks that this subset gives rise to
$|\sigma'(c_j)|$ consecutive letters all equal to $\theta(x)$ in
$w_\sigma$ and that this is a maximal set of consecutive repeated
letters. This shows that one can recover $\sigma'$.  To see that $w$
is recovered correctly, one notes that if $\sigma'(c_j)$ and
$\sigma'(c_{j+1})$ contain the same letter $r$ then $w_j < w_{j+1}$
so by definition $w_j$ is placed correctly before $w_{j+1}$ in
$w^{(r)}_\sigma$.

\end{proof}


\begin{example} We illustrate the proof of Theorem~\ref{thm:multiPpart}.  Let $\theta_s$ be the labeling
$$
\tableau{{3}&{4}&{5} \\ {1}&{2}}
$$
of the shape $\lambda = (3,2)$ as in Example \ref{ex:s}.  Take the
$(\lambda,\theta_s)$-partition \setcellsize{20}
$$
\tableau{{12}&{235}&{5678} \\ {45}&{8}}
$$
in $\AA(\lambda,\theta_s)$.  Then we have $$w_{\sigma} =
(3;3,4;4;1;1,4,5;5;5;2,5),$$ where for example
$w_{\sigma}^{(2)}=(3,4)$ since the cells labeled $3$ and $4$ contain
the number $2$ in $\sigma$. Therefore $$w = (3,4,1,4,5,2,5)$$ and
the corresponding composition $\C(w)$ is $(2,3,2)$.  Then $\sigma'$
if written as sequence is
$$\{1,2\},\{2,3\},\{4,5\},\{5\},\{5,6,7\},\{8\},\{8\}.$$ For example
$\sigma'(c_1) = \{1,2\}$ since $w_{\sigma}^{(1)}$ and
$w_{\sigma}^{(2)}$ contribute $3$'s into the beginning of
$w_{\sigma}$.

The inverse map $\beta$ can now be understood as follows: we parse
$w$ and $\sigma'$ in parallel and place $\sigma'(c_i)$ into the cell
$\theta_s^{-1}(w_i)$.  For example we place $\{1,2\}$ into the cell
labeled $3$, $\{2,3\}$ into the cell labeled $4$, and so on.

\end{example}

\begin{example} \label{ex:mJH}
We give an example of the decomposition of $\K_{\lambda,\theta_s}$
into multi-fundamental quasisymmetric functions. Take $\lambda =
(3,1)$ and take the labeling $\theta_s$ of the cells as described in
Example \ref{ex:s}:
\setcellsize{15}
$$
\tableau{{2}&{3}&{4}\\{1}}
$$

%

Then the usual Jordan-Holder set consists of the sequences
$(2,1,3,4)$, $(2,3,1,4)$, and $(2,3,4,1)$. This coincides with the
part $\tilde \J_4(\lambda,\theta_s)$ of the multi-Jordan-Holder set.
The set $\tilde \J_5(\lambda,\theta_s)$ consists of the words
$(2,1,3,1,4)$, $(2,1,3,4,1)$, $(2,3,1,3,4)$, $(2,3,1,4,1)$, and
$(2,3,4,1,4)$. This gives us the following part of the decomposition
$$\K_{(3,1),\theta_s} =
\L_{(1,3)}+\L_{(2,2)}+\L_{(3,1)}+\L_{(1,2,2)}+\L_{(1,3,1)}+\L_{(2,3)}+\L_{(2,2,1)}+\L_{(3,2)}
+ \cdots$$
\end{example}

\subsection{Structure of $\mQSym$}
Let $\alpha$ be a composition of $n$.  We let $w(\alpha)$ denote any
permutation such that $\alpha = \C(w(\alpha)):=\C(\Des(w(\alpha)))$.

\begin{prop} \label{prop:Lshuffle}
Let $\alpha$ be a composition of $n$ and $\beta$ be a composition of
$m$. Then
$$
\L_\alpha \, \L_\beta = \sum_{u \in {\rm Sh}^\m(w(\alpha),w(\beta))}
\L_{\C(u)},
$$
where ${\rm Sh}^\m(w(\alpha),w(\beta))$ denotes the set of
multishuffles of $w(\alpha)$ and $w(\beta)+n:= (w_1(\beta) +
n)(w_2(\beta) + n) \cdots (w_m(\beta) + n)$.
\end{prop}

\begin{proof}
Take two chains $C = c_1 < c_2 < \ldots < c_n$ and $C' = c'_1 < c'_2
< \ldots < c'_m$.  We label the disjoint union poset $C \cup C'$ by
setting $\theta(c_i) = w_i(\alpha)$ and $\theta(c'_i) = w_i(\beta) +
n$.  Then the multi-Jordan-Holder set $\tilde \J(C\cup C',\theta)$
is exactly the set of multishuffles of $w(\alpha)$ and $w(\beta)+n$.
Since $\K_{C,w(\alpha)} \, \K_{C',w(\beta)} = \K_{C \cup C',\theta}$
we obtain the claimed result by Theorem~\ref{thm:multiPpart}.
\end{proof}

If $f(x)$ is a formal linear combination of the multi-fundamental
quasisymmetric functions $\L_\alpha(x)$, we let $f(x,y)$ denote the
corresponding formal power series in the variables $x_1,
x_2,\ldots,y_1,y_2,\ldots$ obtained by considering $f(x)$ as a
quasisymmetric function (of unbounded degree).

\begin{prop}\label{prop:Lcoprod}
Let $\alpha$ be a composition.  Then
$$
\L_\alpha(x,y) = \sum_{u \otimes u' \in {\rm Cuut}^\m(w(\alpha))}
\L_{\C(u)}(x) \otimes \L_{\C(u')}(y)
$$
where ${\rm Cuut}^\m(w(\alpha))$ denote the set of terms in the cuut
coproduct of $w(\alpha)$.
\end{prop}
\begin{proof}
The power series $\L_\alpha(x,y)$ is the weight generating function
of set-valued $P$-partitions of a labeled chain $(C,w(\alpha))$
using the ordered set of letters $1 < 2 < \cdots < 1' < 2' < \cdots$
where the unprimed letters $i$ are given weight $x_i$ and the primed
letters $i'$ are given weight $y_i$.  There are two kinds of such
set-valued $P$-partitions $\theta$: either (1) there is some $k \in
[1,|\alpha|]$ so that $\theta(c_i)$ contains only unprimed letters
for $i \leq k$ and $\theta(c_j)$ contains only primed letters for $j
> k$, or (2) there is a unique $k \in [1,|\alpha|]$ so that
$\theta(c_k)$ uses both primed and unprimed letters.  This gives
rise to the two kinds of terms in the cuut coproduct: (1)
corresponds to the terms which occur in the usual cut coproduct
while (2) corresponds to the extra terms obtained by cutting in the
middle of $w(\alpha)_k$.
\end{proof}

\begin{theorem}\label{thm:mQSym}The map $\psi:\mMPR \to \prod_\alpha \Z \L_\alpha$
given by $w \mapsto L_{\C(w)}$ is a bialgebra morphism, identifying
$\prod_\alpha \Z \L_\alpha$ with $\mQSym$.
\end{theorem}
\begin{proof}
It is clear that the kernel of $\psi$ is the ideal $I$ of
Lemma~\ref{lem:mQSym}.  The fact that $\psi$ is a bialgebra morphism
follows from Propositions \ref{prop:Lshuffle} and \ref{prop:Lcoprod}
and the proof of Lemma~\ref{lem:mQSym}.
\end{proof}

\subsection{Further properties of the multi-fundamental quasisymmetric functions}
We begin by describing how to express $\L_{\alpha}$ as a (infinite)
linear combination of $L_\alpha$'s.  For a formal power series
$f(x_1,x_2,\ldots)$ of possibly unbounded degree we let $H_i(f)$
denote the homogeneous component of degree $i$.  We define a family
of linear maps $^{(i)}: \QSym \to \QSym$ by ${L^{(i)}_{\alpha}} =
H_{i+|\alpha|}(\L_{\alpha})$, and extending by linearity.  In
particular $f^{(0)}=f$ for any $f \in \QSym$.

Suppose $D \subset [n-1]$ is a subset thought of as a descent set
and $E \subset [n+i-1]$.  An injective and order-preserving map
$t:[n-1] \to [n+i - 1]$ is an {\it $i$-extension} of $D$ to $E$ if
$t(D) \subset E$ and $(E \backslash t(D)) = ([n+i - 1] \backslash
t([n-1])$.  In other words, $E$ is the union of $t(D)$ and the set
of elements not in the image of $t$ (in particular $E$ is determined
by $t$ and $D$). An immediate consequence is that $|E| = |D| + i$.
Note that there may be many $i$-extensions even when $D$ and $E$ are
fixed. For example, if $D = [2] \subset [2]$ and $E = [3] \subset
[3]$ then there are three $1$-extensions of $D$ to $E$,
corresponding to the three injective order-preserving maps $t:[2]
\to [3]$.  We denote the set of $i$-extensions of $D$ to $E$ by
$T(D,E)$.

%
%


\begin{theorem}\label{thm:pump}
Let $\alpha$ be a composition $n$ and $D = D(\alpha)$ the
corresponding descent set.  Then for each $i \geq 0$, we have
\begin{equation}\label{eq:L}
L^{(i)}_{\alpha} = \sum_{E \subset [n+i - 1]} |T(D,E)|\; L_{\C(E)}
\end{equation}
and
\begin{equation}\label{eq:M}
M^{(i)}_{\alpha} = \sum_{E \subset [n+i - 1]} |T(D,E)|\; M_{\C(E)}.
\end{equation}
For each $i, j \geq 0$ and $f \in \QSym$, one has $(f^{(i)})^{(j)} =
{i+j \choose i} f^{(i+j)}$.
\end{theorem}

\begin{proof}
Let $w = w(\alpha)$ and consider the subset $\AA_i(C,w) \subset
\AA(C,w)$ consisting of set-valued $(C,w)$-partitions $\sigma$ of
size $|\sigma| = n + i$.  We must show that the generating function
of $\AA_i(C,w)$ is equal to $\sum_{E \subset [n+i - 1]} |T(D,E)|\;
L_{\C(E)}$.  Indeed for each pair $t \in T(D,E)$ for some $E$, the
function $L_{\C(E)}$ is the generating function of all $\sigma \in
\AA_i(C,w)$ satisfying $|\sigma(c_i)| = t(i) - t(i-1)$ where one
defines $t(0) = 0$ and $t(n) = n+i$.  Indeed one obtains a (usual)
$(C',w(\C(E))$-partition $\sigma' \in {\mathcal A}(C',w(E))$ by
assigning the elements of $\sigma(c_i)$ in increasing order to
$c'_{t(i-1)+1}, \ldots,c'_{t(i)}$, where $C' = c'_1 < c'_2 < \ldots
< c'_{n+i}$ is a chain with $n + i$ elements. This proves
(\ref{eq:L}).

To prove (\ref{eq:M}) it suffices to show that (\ref{eq:M}) implies
(\ref{eq:L}) since both $\{L_\alpha\}$ and $\{M_\alpha\}$ form bases
of $\QSym$.  Assuming (\ref{eq:M}), we calculate
$$
L_{\C(D)}^{(i)} = \sum_{C \supset D}M_{\C(C)}^{(i)} = \sum_{C
\supset D} \sum_{B \subset [n+i-1]} |T(C,B)|\;M_{\C(B)}.
$$
We show that $$\sum_{C \supset D} \sum_{B \subset [n+i-1]}
|T(C,B)|\;M_{\C(B)} = \sum_{E \subset[n+i-1]} |T(D,E)| \sum_{B
\supset E}M_{\C(B)}$$ from which our claim will follow.  A term
$M_{\C(B)}$ on the right hand side is indexed by the following data:
an $i$-extension $t:[n-1] \to [n+i-1]$ (of $D$) and the set $B
\backslash E$ contained in $t([n-1] \backslash D)$.  But since $t$
is injective, this is the same as giving the subset $t^{-1}(B
\backslash E) \subset [n-1]\backslash D$ and the $i$-extension
$t:[n-1] \to [n+i-1]$.  This is exactly the information indexing
terms on the left hand side.

To prove the last claim, let $F \subset [n+i+j-1]$ and $t:[n-1] \to
[n+i+j-1]$ be an $i+j$-extension of $D$ to $F$.  Now let $[n+i+j-1]
\backslash t([n-1]) = S \cup S'$ be a decomposition into a set $S$
containing $j$-elements and $S'$ containing $i$-elements. Then there
is a unique $E \subset[n+i-1]$ and $t':[n+i-1] \to [n+i+j-1]$ such
that $t'$ is a $j$-extension of $E$ to $F$ satisfying $S = [n+i+j-1]
\backslash t'([n+i-1])$. Furthermore, there is a unique $t'': [n-1]
\to [n+i-1]$ which is an $i$-extension of $D$ to $E$ such that
$[n+i-1] \backslash t''([n-1]) = (t')^{-1}(S')$.  The composition
$t' \circ t''$ is equal to $t$.  The correspondence $(t,S) \mapsto
(t',t'')$ is a bijection. Since there are exactly ${i+j \choose i}$
choices for $S$, this proves that $(f^{(i)})^{(j)} = {i+j \choose i}
f^{(i+j)}$.
\end{proof}

\begin{example}
Take $\alpha = (2,1)$. Then $D(\alpha)=\{2\}$ and we have the
following equality: $$\L_{\alpha}^{(2)} = L_{(1,1,2,1)}+2
L_{(1,2,1,1)} + 3 L_{(2,1,1,1)}.$$ Here for example $|T(\{2\},
\{1,3,4\})|=2$ since there are two maps $t_1: \{1,2\} \mapsto
\{2,3\}$ and $t_2: \{1,2\} \mapsto \{2,4\}$ which satisfy the needed
condition. By applying the rule we get the following part of
decomposition of $\L_{(2,1)}$ into $L$-s: $$\L_{(2,1)} =
L_{(2,1)}+L_{(1,2,1)}+2 L_{(2,1,1)}+L_{(1,1,2,1)}+2 L_{(1,2,1,1)} +
3 L_{(2,1,1,1)} + \ldots.$$
\end{example}

\begin{remark}
One may define the {\it {multi-monomial quasisymmetric functions}}
$\tilde M_\alpha$ by analogy with usual monomial symmetric functions
$M_\alpha$:
$$\tilde M_\alpha = \sum_{D(\beta) \subset
D(\alpha)}(-1)^{|D(\alpha)|-|D(\beta)|} \L_{\beta}$$ where the
summation is restricted to compositions satisfying $|\beta| =
|\alpha|$.  It is clear that $H_{|\alpha|}(\tilde
M_{\alpha})=M_{\alpha}$.
\end{remark}

As an application of Theorem~\ref{thm:pump}, we give a curious
property of descents in the multi-Jordan-Holder set for partition
shapes.  It generalizes the following statement concerning usual
Jordan-Holder sets.

\begin{theorem} \cite[Theorem 7.19.9]{EC2} \label{thm:StaDes}
Let $\lambda/\mu$ be a skew Young diagram with $n$ boxes.  For any
$1 \leq i \leq n-1$ the number $d_i(\lambda/\mu)$ of $w \in
\J(\lambda/\mu,\theta_{s})$ for which $i \in D(w)$ is independent of
$i$.
\end{theorem}

Theorem~\ref{thm:StaDes} is proved essentially using the following
lemma, implicit in the argument of \cite{EC2}.  Define the linear
transformation $\phi_i: \QSym \longrightarrow \mathbb \Z$ by

$$
\phi_i(L_{\alpha}) =
\begin{cases}
1 & \text{if $i \in D(\alpha)$}\\
0 & \text{otherwise,}
\end{cases}
$$
and extending by linearity.


\begin{lemma} \label{lem:balanced}
Let $f = \sum c_{\alpha} M_{\alpha} \in \QSym$ be a homogeneous
quasisymmetric function of degree $n$.  Then $\phi_i(f)$ does not
depend on $i$ if and only if $c_{(2,1,\ldots, 1)} = c_{(1,2,\ldots,
1)} = \ldots = c_{(1,1,\ldots, 2)}$.
\end{lemma}
\begin{proof}
The statement follows from the equation $L_\alpha = \sum_{E \supset
D(\alpha)} M_{\C(E)}$.
\end{proof}

We call $f$ satisfying the condition of Lemma~\ref{lem:balanced}
{\it {balanced}}.  We also call the monomials labeled by the
compositions $\alpha = (1^i,2,1^{n-2-i})$ {\it {balancing}}.

\begin{lemma}\label{lem:pumpbalance}
If $f$ is a homogeneous balanced quasisymmetric function, then so is
$f^{(i)}$ for any $i \geq 0$.
\end{lemma}

\begin{proof}
The claim clearly holds for $i = 0$.  Using the last statement of
Theorem~\ref{thm:pump}, we may assume that $i = 1$.  Suppose $f$ has
degree $n$.  Let $D_j = [n-1] \backslash\{j\}$ for $1 \leq j \leq n
- 1$ and $E_k = [n] \backslash \{k\}$ for $1 \leq k \leq n$.  By
Theorem~\ref{thm:pump}, to calculate the coefficient of
$M_{\C(E_k)}$ in $f^{(1)}$ it suffices to find the $1$-extensions
$t:[n-i] \to [n]$ of some $D \subset [n-1]$ to $E_k$. But $|E_k| =
n+1$ so such $D$ satisfy $|D| = n - 2$ and so must be of the form $D
= D_j$ for some $j$.  The number of $1$-extensions of $D_j$ to $E_k$
is equal to $n - k$ if $j = k$, equal to $k - 1$ if $j = k-1$, and
equal to 0 otherwise.  By Lemma~\ref{lem:balanced} the coefficient
of $M_{\C(D_j)}$ in $f$ does not depend on $j$, thus the coefficient
of $M_{\C(E_k)}$ in $f^{(1)}$ is equal to $n-1$ times the
coefficient of $M_{\C(D_j)}$ in $f$, which does not depend on $k$.
Again by Lemma~\ref{lem:balanced}, $f^{(1)}$ must be balanced.
\end{proof}

We prove the following generalization of the Theorem
\ref{thm:StaDes}.

\begin{theorem} \label{thm:m-J-H}
Let $|\lambda/\mu| = n$. For any $N \geq n$ and $1 \leq i \leq N-1$
the number $d_i(\lambda/\mu)$ of $w \in \tilde
\J_N(\lambda/\mu,\theta_{s})$ for which $i \in D(w)$ is independent
of $i$.
\end{theorem}

\begin{proof}
Let $G'_{\lambda/\mu} = \K_{\lambda/\mu,\theta_s}$, which has lowest
degree homogeneous component $H_0(G'_{\lambda/\mu})=
s_{\lambda/\mu}$ equal to a Schur function.  We shall see later
(Corollary~\ref{cor:Gsymmetric})) that $G'_{\lambda/\mu}$ is a
symmetric function (of unbounded degree). By Theorem
\ref{thm:multiPpart} for each $N \geq n$ we have

$$
H_N(G'_{\lambda/\mu}) = \sum_{n \leq m \leq N} \sum_{w \in \tilde
\J_m(\lambda/\mu,\theta_{s})} L_{\C(w)}^{(N-m)}.
$$

We know that $s_{\lambda/\mu} = \sum_{w \in \tilde
\J_n(\lambda/\mu,\theta_{s})} {L_{\C(w)}}$  is symmetric, and thus
balanced by Lemma~\ref{lem:balanced}. By
Lemma~\ref{lem:pumpbalance}, $s_{\lambda/\mu}^{(1)}$ is balanced and
we also know $H_{n+1}(G'_{\lambda/\mu})$ is symmetric so $\sum_{w
\in \tilde \J_{n+1}(\lambda/\mu,\theta_{s})} {L_{\C(w)}} =
H_{n+1}(G'_{\lambda/\mu}) - s_{\lambda/\mu}^{(1)}$ is also balanced.
Proceeding in this manner we conclude that each of the sums $\sum_{w
\in \tilde \J_{N}(\lambda/\mu,\theta_{s})} {L_{\C(w)}}$ is balanced.
By Lemma \ref{lem:balanced} we obtain exactly the needed result.
\end{proof}

\begin{example}
One can check that in Example \ref{ex:mJH} for each $1 \leq i \leq
4$ there exists exactly two elements of $\tilde
\J_5(\lambda,\theta_s)$ with $i$ as a descent.
\end{example}

\section{$K$-theory of Grassmannians and $\mSym$}
\label{sec:mSym}

\subsection{Fomin-Greene operators}
Let $\Lambda$ denote the set of partitions as before.  If $\lambda =
(\lambda_1 \geq \lambda_2 \geq \cdots \geq \lambda_l > 0)$ is a
partition, then it contains the boxes $(i,j)$ for $1 \leq i \leq l$
and $1 \leq j \leq \lambda_i$.  The box $(i,j)$ is on {\it diagonal}
$j - i$.  We say that $\lambda$ has an {\it inner corner} on the
$i$-th diagonal if there exists $\mu \in \Lambda$ such that
$\lambda/\mu$ is a single box on the $i$-th diagonal.  Similarly,
$\lambda$ has an {\it outer corner} on the $i$-th diagonal if there
exists $\mu \in \Lambda$ such that $\mu/\lambda$ is a single box on
the $i$-th diagonal.

Fix a partition $\nu \in \Lambda$.  Let $\Z\Lambda_\nu = \oplus_{\nu
\subset \lambda} \Z\cdot \lambda$ denote the free $\Z$-module with a
basis of partitions containing $\nu$, equipped with a non-degenerate
pairing $\ip{.,.}: \Z\Lambda_\nu \times \Z\Lambda_\nu \to \Z$
defined by $\ip{\lambda,\mu} = \delta_{\lambda\mu}$. Now for each $i
\in \Z$, define a $\Z$-linear operator $v_i:\Z\Lambda_\nu \to
\Z\Lambda_\nu$ by
$$
v_i \cdot \lambda = \begin{cases} \mu & \mbox{if $\lambda$ has an
outer corner $\mu/\lambda$ on the $i$-th diagonal,} \\
\lambda & \mbox{if $\lambda$ has an inner corner not contained in
$\nu$ on the $i$-th
diagonal,} \\
0 & \mbox{otherwise,}
\end{cases}
$$
and extending by linearity.  The operators $v_i$ satisfy the
relations:
\begin{align*}
v_i^2 &= v_i & \mbox{for each $i \in \Z$,} \\
v_iv_{i+1} v_i &= v_{i+1}v_iv_{i+1} = 0 & \mbox{for each $i \in
\Z$,}
\\
v_i v_j &= v_j v_i &\mbox{for each $i,j \in \Z$ with $|i-j| \geq
2$.}
\end{align*}

Operators very closely related to the $v_i$ are studied by Fomin and
Greene~\cite{FG} and operators differing from ours by a sign are
studied by  Buch~\cite{B}. We briefly explain their connection
with set-valued $P$-partitions to draw an analogy with our
construction of $\MSym$ later.  Define a formal power series
$$
A(x) = \cdots (1+xv_2)(1+xv_1)(1+xv_0)(1+xv_{-1}) \cdots
$$
with coefficients in operators on $\Z\Lambda_\nu$.  The action of
$A(x)$ on $a \in \Z\Lambda_\nu$ gives a well defined element of
$\Z\Lambda_\nu[[x]]$. The following result is essentially
\cite[Theorem 3.1]{B} with a sign omitted.
\begin{lemma}\label{lem:K}
Let $\lambda/\mu$ be a skew partition.  Then
$$
\K_{\lambda/\nu,\theta_s}(x_1,x_2,\ldots) = \ip{\cdots A(x_2)A(x_1)
\cdot \nu,\lambda}.
$$
\end{lemma}

The next result follows from work of Fomin and Greene~\cite{FG}.
\begin{lemma}\label{lem:commute}
We have $A(x)A(y) = A(y)A(x)$ as operators on $\Z\Lambda_\nu$.
\end{lemma}

\begin{corollary}\label{cor:Gsymmetric} Let $\lambda/\mu$ be a skew shape.  Then
$\K_{\lambda/\mu,\theta_s}$ is a symmetric function (of unbounded
degree).
\end{corollary}
\begin{proof}
Follows immediately from Lemmas~\ref{lem:K} and \ref{lem:commute}.
\end{proof}

\subsection{The Hopf algebra $\mSym$ and $K$-theory of Grassmannians}
Let $\mSym = \prod_\lambda \Z \cdot \K_{\lambda,\theta_s}$ denote
the subspace of $\mQSym$ continuously spanned by the generating
functions $\K_{\lambda,\theta_s}$ as $\lambda$ varies over all
partitions.  For a fixed composition $\alpha$, $\L_\alpha$ only
occurs in finitely many $\K_{\lambda,\theta_s}$ so $\mSym$ is indeed
a subspace of $\mQSym$.  For convenience we now write $\K_\lambda$
instead of $\K_{\lambda,\theta_s}$.  Also we shall call a set-valued
$(\lambda/\mu,\theta_s)$-partition $\sigma$ simply a {\it set-valued
tableau} of shape $\lambda/\mu$.

\begin{prop}\label{prop:mSymSym}
The space $\mSym$ is a Hopf subalgebra of $\mQSym$.  It is
isomorphic to the completion $\widehat{\Sym} = \prod_\lambda \Z
\cdot s_\lambda$ of the algebra of symmetric functions.
\end{prop}
\begin{proof}
Since the lowest degree homogeneous component of $\K_{\lambda}$ is
equal to the Schur function $s_\lambda$, the space $\mSym$ is equal
to the space $\prod_\lambda \Z \cdot s_\lambda$ of arbitrary
linearly combinations of Schur functions. Thus $\mSym$ is the Hopf
subalgebra of $\mQSym$ consisting of symmetric functions of
unbounded degree.
\end{proof}

In~\cite{B}, Buch studied a bialgebra $\Gamma = \oplus_\lambda \Z
\cdot G_{\lambda}$ spanned by the set $\{G_\lambda\}$ of stable
Grothendieck polynomials.  The stable Grothendieck polynomials were
first studied by Fomin and Kirillov \cite{FK} and defined as stable
limits of Grothendieck polynomials \cite{LS}. For our purposes, they
can be defined as follows.

\begin{theorem} \cite[Theorem 3.1]{B} \label{thm:sv}
The stable Grothendieck polynomial $G_{\nu/\lambda}(x)$ is given by
the formula $$G_{\nu/\lambda}(x) = \sum_T
(-1)^{|T|-|\nu/\lambda|}x^T,$$ where the sum is taken over all
set-valued tableaux $T$ of shape $\nu/\lambda$.
\end{theorem}

Thus the symmetric function $G_\lambda$ can be obtained from
$\K_\lambda$ by changing the degree $n$ homogeneous component of
$\K_\lambda$ by the sign $(-1)^{n-|\lambda|}$, or in other words one
has $\K_\lambda(x_1,x_2,\ldots) =
(-1)^{|\lambda|}\,G_\lambda(-x_1,-x_2,\ldots)$.

Buch related the structure constants of $\Gamma$ to the $K$-theory
$K^\circ\Gr(k,\c^n)$ of the Grassmannian $\Gr(k,\c^n)$ of $k$-planes
in ${\mathbb C}^n$. In addition, Buch described the structure
constants completely using the combinatorics of set-valued tableaux.
We briefly describe the connections with $K$-theory here and return
to the combinatorial descriptions from a dual point of view later
(see Section~\ref{sec:khom}).

Let $R = (n-k)^k$ denote the the rectangle with $k$ rows and $n-k$
columns and let $I_R$ be the subspace of $\Gamma$ spanned by
$G_\lambda$ for $\lambda$ not contained in $R$.  Let
$K^\circ\Gr(k,\c^n)$ denote the Grothendieck group of algebraic
vector bundles on $\Gr(k,\c^n)$.  It is naturally isomorphic to the
Grothendieck group $K_\circ\Gr(k,\c^n)$ of coherent sheaves on
$\Gr(k,\c^n)$.  The $K$-group $K^\circ\Gr(k,\c^n)$ is spanned by the
classes $[\O_\lambda]$ of structure sheaves of Schubert varieties
$X_\lambda$ indexed by partitions $\lambda \subset R$.  For
convenience we set $[\O_\lambda] = 0$ if $\lambda$ does not fit
inside $R$.  The $K$-theory $K^\circ\Gr(k,\c^n)$ becomes a
commutative ring when equipped with the multiplication induced by
tensor products of vector bundles.

\begin{theorem}[{\cite[Theorem 8.1]{B}}] \label{thm:Buch}
The map $G_\lambda \mapsto [\O_\lambda]$ induces an isomorphism of
rings $\Gamma/I_R \simeq K^\circ\Gr(k,\c^n)$.
\end{theorem}

Note that $\Gamma/I_R$ is isomorphic to the quotient of $\mSym$ by
the continuous span $\prod_{\lambda \nsubseteq R}\K_\lambda$ since
both spaces are finite-dimensional.  To explain the geometric
meaning of the coproduct, fix $k_1 < n_1$ and $k_2 < n_2$.  Taking
the direct sum of vector spaces induces a map $\phi:
\Gr(k_1,\c^{n_1}) \times \Gr(k_2,\c^{n_2}) \to
\Gr(k_1+k_2,\c^{n_1+n_2})$.  Then Buch shows that
$\phi^*([\O_\lambda]) = \sum_{\mu,\nu} d^{\mu\nu}_\lambda [\O_\mu]
\otimes [\O_\nu]$ where we have identified $K^\circ
\Gr(k_1,\c^{n_1}) \times \Gr(k_2,\c^{n_2})$ with
$K^\circ\Gr(k_1,\c^{n_1}) \otimes K^\circ\Gr(k_2,\c^{n_2})$.  Here
the $d^{\mu\nu}_\lambda$ are the structure constants of the
coproduct: $\Delta G_\lambda = \sum_{\mu,\nu} d^{\mu\nu}_\lambda
G_\mu \otimes G_\nu$.

\section{The big Multi Malvenuto-Reutenauer Hopf algebra} \label{sec:MMPR}
\subsection{Big Multi-permutations and set
compositions}\label{sec:Mperm}
\begin{definition}
A {\it big Multi-permutation} or $\M$-permutation of $[n]$ is a word
$w = w_1 w_2 \cdots w_k$ such that (a) each $w_i$ is a subset of
$\P$ not containing consecutive numbers and (b) the disjoint union
$\sqcup_{i = 1}^k w_i$ is equal to the set $[n]$. We say that $w$
has length $\ell(w) = k$.
\end{definition}

Denote the set of $\M$-permutations of $[n]$ by $\bS_{n}$ and let
$\bS_\infty = \cup_{n} \bS_{n}$.  For example, we have $w =
[(1,3),(5,7,9),(10),(4,6),2,8] \in \bS_{10}$ which has length
$\ell(w) = 6$.  By convention $\bS_0$ contains a single element --
the empty $\M$-permutation $\emptyset$.  We let $\MMPR = \oplus_{w
\in \bS} \Z.w$ denote the free $\Z$-module of finite $\Z$-linear
combinations of $\M$-permutations.

Recall that a {\it set composition} $w = w_1 w_2 \cdots w_k$ of a
finite set $S$ is a word $w = w_1 w_2 \cdots w_k$ such that the
disjoint union $\sqcup_{i = 1}^k w_i$ is equal to the set $S$.  Thus
$w_i$ may contain consecutive numbers.  If $w = w_1 w_2 \cdots w_k$
is a set composition of $S \subset \P$, we define the {\it
standardization} $\st(w) \in \bS_\infty$ by repeatedly doing the
following operations until one has a $\M$-permutation:
\begin{enumerate} \item[R1] delete the letter $i+1$ if both $i$ and $i+1$ belong
to some $w_j$, and
\item[R2] reduce all letters larger than $i$ by 1
if $i$ is not present in any $w_j$.
\end{enumerate}
It is clear that $\st(w)$ is a well defined $\M$-permutation. Also
if $w$ is a set composition of $S$ and $T \subset S$ then the {\it
restriction} $w|_T$ of $w$ to $T$ is obtained by intersecting each
$w_i$ with $T$ and removing all the $w_j$ which become empty. The
restriction $w|_T$ is a set composition of $T$.

Let $w = w_1 \cdots w_k \in \bS_m$ and $u = u_1 \cdots u_l \in
\bS_n$.  Define a product $\bullet$ on $\MMPR$ by extending the
formula
$$w \bullet u =
\sum v$$ by linearity, where the (finite) sum is taken over
\begin{enumerate}
\item
all $v \in \bS_{m+n}$ such that $v|_{[m]} = w$ and
$\st(v|_{[m+1,m+n]}) = u$; and
\item
all $v \in \bS_{m+n-1}$ such that $v|_{[m]} = w$ and
$\st(v|_{[m,m+n-1]}) = u$.
\end{enumerate}

\begin{prop}\label{prop:MMPRassociative}
The product $\bullet: \MMPR \otimes \MMPR \to \MMPR$ is associative.
\end{prop}
\begin{proof}
Let $w = w_1 \cdots w_k \in \bS_m$, $u = u_1 \cdots u_l \in \bS_n$,
and $x = x_1 \cdots x_r \in \bS_p$.  Then one checks that both $(w
\bullet u) \bullet x$ and $w \bullet (u \bullet x)$ are equal to the
sum over
\begin{enumerate}
\item
all $v \in \bS_{m+n+p}$ such that $v|_{[m]} = w$,
$\st(v|_{[m+1,m+n]}) = u$ and \\ $\st(v|_{[m+n+1,m+n+p]}) = x$; and
\item
all $v \in \bS_{m+n+p-1}$ such that $v|_{[m]} = w$,
$\st(v|_{[m,m+n-1]}) = u$ and \\ $\st(v|_{[m+n,m+n+p-1]}) = x$; and
\item
all $v \in \bS_{m+n+p-1}$ such that $v|_{[m]} = w$,
$\st(v|_{[m+1,m+n]}) = u$ and \\ $\st(v|_{[m+n,m+n+p-1]}) = x$; and
\item
all $v \in \bS_{m+n+p-2}$ such that $v|_{[m]} = w$,
$\st(v|_{[m,m+n-1]}) = u$ and \\ $\st(v|_{[m+n-1,m+n+p-2]}) = x$.
\end{enumerate}
If $v \in \bS_{m+n+p-1}$ satisfies both the conditions in (2) and
(3) then it occurs with multiplicity two in the product of $w,u$ and
$x$.
\end{proof}

We can also give an alternative recursive definition of the product
$\bullet$. For two set compositions $a,b$ let $[a,b]$ denote their
concatenation (assuming that the result is a set composition).  We
can extend $[.,.]$ by linearity by distributing it over addition.
We first define the semishuffle product $\circ$ as follows. Let $u =
u_1 u'$ and $v = v_1 v'$ be $\M$-permutations where $u_1,v_1$ are sets and $u',v'$ are
set compositions.  Then $$u \circ v =
[u_1, u' \circ v] + [v_1, u \circ v'] + [u_1 \cup v_1, u'
\circ v'].$$ For a $\M$-permutation $v$, let $(v+n)$ denote the set composition obtained by
increasing every number in $v$ by $n$. Let $u \in S^{\M}_n$.


\begin{prop}
We have $u \bullet v = \st(u \circ (v+n))$.
\end{prop}

Now define the coproduct $\vartriangle$ on $\MMPR$ by
$$\vartriangle w = \sum_{[u,v]=w} \st(u) \otimes \st(v),$$
where the sum is over all pairs $(u,v)$ of (possibly empty) set
compositions which concatenate to $w$.  We extend $\vartriangle$ by
linearity to give $\vartriangle: \MMPR \to \MMPR \otimes \MMPR$. The
following result is immediate from the definition.

\begin{prop}\label{prop:MMPRcoassociative}
The coproduct $\vartriangle: \MMPR \otimes \MMPR \to \MMPR$ is
coassociative.
\end{prop}

Define the unit map $\eta: \Z \to \MMPR$ by $\eta(1) = \emptyset$
and the counit map $\varepsilon: \MMPR \to \Z$ by taking the
coefficient of $\emptyset$.

\begin{theorem}\label{thm:MMPR}
The space $\MMPR$ is a bialgebra with product $\bullet$, coproduct
$\vartriangle$, unit $\eta$, and counit $\varepsilon$.
\end{theorem}
\begin{proof}
It is easy to check that $\MMPR$ is a unital associative algebra and
a counital coassociative coalgebra.  We must therefore check that
$\bullet$ and $\vartriangle$ are compatible.  Let $w = w_1 \cdots
w_k \in \bS_m$ and $u = u_1 \cdots u_l \in \bS_n$.  Then both
$\vartriangle(w \bullet u)$ and $(\vartriangle w) \bullet
(\vartriangle u)$ are sums over terms $a \otimes b$ which can be
described as follows.  First let $u'$ be the unique set composition
of $[m+1,m+n]$ such that $\st(u') = u$.  Then $a = \st(a' = a'_1
a'_2 \cdots a'_r)$ where for each $1 \leq i \leq r$, we have $a'_i$
is either (a) equal to some $w_j$, or (b) equal to some $u'_s$, or
(c) the union of $w_j$ and $u'_s$.  Also if $w_{j_1}$ appears in
$a'_{i_1}$ and $w_{j_2}$ in $a'_{i_2}$ then $i_1 < i_2 \Rightarrow
j_1 < j_2$ and furthermore all $w_j$ with $j < j_2$ have to appear
in $a'$ (and the analogous statement for $u'$).  Similarly $b =
\st(b')$ for an analogously described $b'$, and in addition the
disjoint union of all the numbers in $a'$ and $b'$ is equal to
$[m+n]$.
\end{proof}

Call an element $w \in \bS_\infty$ {\it {irreducible}} if it cannot
be written in the form $w = v / u = v_1 \ldots v_k (u_1 +n) \ldots
(u_l+n)$ for two smaller (that is, non-empty) $\M$-permutations $v,u
\in \bS_\infty$. Let ${\rm Irr}^{\M}$ be the set of irreducible
$\M$-permutations.

\begin{theorem} \label{thm:MMRfree}
The algebra $\bS_\infty$ is free over the set of its irreducible
elements.
\end{theorem}

\begin{proof}
Any $\M$-permutation $w \in \bS_\infty$ can be uniquely written as
$w^1 / w^2 / \cdots / w^k$ where the $w^i \in {\rm Irr}^\M$ are
irreducible; we say that $w$ is $k$-reducible.  Thus the tensor
space $\bigoplus_{n \geq 0}\Z({\rm Irr}^{\M})^{\otimes n}$ is
(naturally) isomorphic as a free $\Z$-module to $\MMPR$.

Now if $w^1,\ldots,w^k \in {\rm Irr}^\M$ are irreducible then it
follows from the definition of $\bullet$ that the only $k$-reducible
term in $w^1 \bullet w^2 \bullet \cdots \bullet w^k$ is $w^1 / w^2 /
\cdots / w^k$. This implies (via a triangularity argument) that the
map $w^1 \otimes w^2 \otimes \cdots \otimes w^k \mapsto w^1 \bullet
w^2 \bullet \cdots \bullet w^k$ induces a surjective algebra
homomorphism $\bigoplus_{n \geq 0}\Z({\rm Irr}^{\M})^{\otimes n} \to
\MMPR$.  The previous paragraph implies that this surjective map is
an isomorphism.

\end{proof}

\begin{remark}
The Hopf structure of $\MMPR$ may be related to the Hopf algebra of
set partitions defined in \cite{PR}.
\end{remark}

\subsection{The antipode of $\MMPR$}
We show that $\MMPR$ has an antipode via a general construction
following \cite[Section 5]{AS} (see also \cite{Tak}). Let $H$ be any
bialgebra with multiplication $m$, coproduct $\Delta$, unit $\mu$
and counit $\varepsilon$.  For each $i \geq 1$ let $m^{(i)}:
H^{\otimes i+1} \to H$ denote the $i$-fold iterated product (by
associativity the order does not matter) and let $\Delta^{(i)}: H
\to H^{\otimes i+1}$ denote the $i$-fold iterated coproduct.  In
addition we set $m^{(0)} = \Delta^{(0)} = {\rm id}: H \to H$,
$m^{(-1)} = \mu$, and $\Delta^{(-1)} = \varepsilon$.  If $f: H \to
H$ is any linear map then its $i$-fold convolution is $f^{(i)} =
m^{(i-1)}f^{\otimes i} \Delta^{(i-1)}$.

Now set $\pi = {\rm id} - \mu \varepsilon$.  If $\pi$ is locally
nilpotent with respect to convolution then the antipode $S: H \to H$
is given by \begin{equation}\label{eq:antipode} S = \sum_{i \geq 0}
(-1)^{i} m^{(i-1)} \pi^{\otimes i} \Delta^{(i-1)}.
\end{equation}

\begin{prop}
The bialgebra $\MMPR$ has an antipode, making it a Hopf-algebra.
\end{prop}
\begin{proof}
By the preceding discussion it suffices to show that $\pi$ is
locally nilpotent. Let $a \in \MMPR$ be non-zero.  Let $n \geq 0$ be
the maximal value such that some $w$ satisfying $\ell(w)= n$ occurs
in $a$ with non-zero coefficient.  Then using the definition of
$\vartriangle$, each term in $\vartriangle^{(n)}a$ must involve
$\emptyset$ in one of its factors. But $\pi(\emptyset) = 0$ so
$\pi^{(n+1)}a = 0$.
\end{proof}

\subsection{Weak order on $\M$-permutations}  The results of this
section suggest that there may be a polytopal structure on
$\M$-permutations.  Let $\SC(n)$ denote the collection of set
compositions of $[n]$ and let $\SC(\infty) = \cup_{n \geq 1}
\SC(n)$.  In \cite{KLNPS} set compositions were considered under the
name of {\it {pseudo permutations}}, while in \cite{PR} they are
identified with faces of permutohedra.  In \cite{KLNPS} and
\cite{PR} independently the following analog of the weak order is
defined  on $\SC(n)$, for each $n \geq 1$.  Let $w = w_1w_2\cdots
w_k \in \SC(n)$.  The covers of $w$ can be completely described in
the following manner.
\begin{enumerate}
\item
if $w_i < w_{i+1}$ then $w \lessdot w_1\cdots w_{i-1} (w_i \cup
w_{i+1}) w_{i+1} \cdots w_k$;
\item
if $w_i = w'_i \sqcup w''_i$ is a disjoint union of non-empty
subsets and $w'_i < w''_i$ then $w \lessdot w_1 \cdots w_{i-1} w''_i
w'_i w_{i+1} \cdots w_k$.
\end{enumerate}
The order $\prec$ on $\SC(n)$ is then the transitive closure of
$\lessdot$.

Now let $\sim$ be the equivalence relation on $\SC(\infty)$ obtained
by taking the transitive closure of the relations $w \sim' \st(w)$
for each $w \in \SC(\infty)$.  Thus each equivalence class contains
a unique $\M$-permutation.  For example,
$$[(1,4,5),7,(2,8,9),(6,10),3] \sim
[(1,4),(6,7,8),(2,9),(5,10,11),3]$$ since
\begin{align*}\st([(1,4,5),7,(2,8,9),(6,10),3]) &=
\st([(1,4),(6,7,8),(2,9),(5,10,11),3])  \\
&=[(1,4),6,(2,7),(5,8),3].\end{align*}

If $w \sim u$ are two equivalent set compositions, let us say that
$w$ {\it contains} $u$ if $u$ can be obtained from $w$ by partially
standardizing: in other words, by using the operations (R1) and (R2)
of Section~\ref{sec:Mperm} a number of times.  The following two
lemmata are easy to establish from the definitions.
\begin{lemma}\label{lem:contain}
If $w \sim u$ are two set compositions then there exists a set
composition $v$ equivalent to both which contains both.
\end{lemma}
\begin{lemma}\label{lem:cover}
Suppose $w \lessdot v$ is a cover relation in $SC(n)$.  If $w'$
contains $w$ then there is a canonical cover relation $w' \lessdot
v'$ such that $v'$ contains $v$.
\end{lemma}

Now define the {\it weak order} $<$ on $\bS_\infty$ by taking the
transitive closure of the following relations: $w \in \bS_\infty$ is
less than $v \in \bS_\infty$ if there exists set compositions $w'$
and $v'$ so that $w' \sim w$, $v' \sim v$ and $w' \lessdot v'$.

\begin{theorem}
The weak order on $\bS_\infty$ is a valid partial order.
\end{theorem}

\begin{proof}
We need to show that for two $\M$-permutations $w,w'$ that if $w <
w'$ and $w' < w$ then $w = w'$.  Alternatively we need to show that
there is no sequence $w_1, w_2, \ldots, w_n, w_1$ of
$\M$-permutations such that $v_1 \lessdot u_2, v_2 \lessdot u_3,
\cdots v_{n-1} \lessdot u_n, v_n \lessdot u_1$ for some set
compositions $v_i \sim w_i$ and $u_i \sim w_i$. Using Lemma
\ref{lem:contain} and Lemma \ref{lem:cover} repeatedly, we may
assume that $v_i = u_i$ for $i \neq 1$.  Thus we are reduced to
proving that if $v < u$ in $\SC(n)$ for some $n$ then $\st(v) \neq
\st(u)$.

Suppose that $v < u$ and $\st(v) = \st(u) = w$ for some
$\M$-permutation $w \in \bS_m$.  Define the ``standardization''
function $f_v:[n] \to [m]$ by requiring that (a) $f_v$ is
increasing, and (b) $f_v^{-1}(i)$ is a non-empty interval completely
contained in the set $v_k$ if $i \in w_k$.  Similarly define $f_u$.
Now since $\st(v) = \st(u) = w$ we must have that $v = v_1 \cdots
v_k$ and $u = u_1 \cdots u_k$ have the same length. Thus it makes
sense to ask for the smallest integer $i$ such that $i \in v_j$ and
$i \in u_k$ for $j \neq k$.  But since $\st(v) = \st(u)$, the letter
$i-1$ must be either in both $v_j$ and $u_j$ or in both $v_k$ and
$u_k$. Let us say that a set composition $x$ has an {\it inversion}
at $(i,j)$ for $i < j$ if $j$ precedes $i$; and a {\it
half-inversion} if $j$ and $i$ belong to the same set of $x$.  We
note that weak order on $\SC(n)$ either increases or does not change
inversions for each $(i,j)$.

Now suppose $j < k$.  If $i-1 \in v_j \cap u_j$ then $v$ has a
half-inversion at $(i-1,i)$ while $u$ has no inversion which is
impossible.  If $i-1 \in v_k \cap u_k$ then $v$ has an inversion at
$(i-1,i)$ while $u$ only has a half-inversion which again is
impossible.  Thus $j > k$.

Suppose first that $i-1 \in v_k \cap u_k$.  Let $i_1$ be the
smallest integer greater than $i$ such that $i_1 \notin u_k$.  Then
$i_1 \in f_u^{-1}(f_v(i))$ so $i_1 \in u_j$.  Since $u$ has no
inversion at $(i,i_1)$, it must also be the case that $v$ has no
inversion at $(i,i_1)$.  Thus $i_1 \in v_{j_1}$ for $j_1 > j$.  Now
let $i_2$ be the smallest integer in $f_u^{-1}(f_v(i_1))$.  Clearly
$i_2 > i_1$ and $i_1 \in u_{j}$ while $i_2 \in u_{j_1}$.  Again $u$
has no inversion at $(i_1,i_2)$ so $v$ has no inversion at
$(i_1,i_2)$ and we deduce that $i_2 \in v_{j_2}$ for $j_2 > j_1$.
Continuing in this manner we produce an infinite sequence
$i_1,i_2,\ldots$.  Since $n$ is finite, we arrive at a
contradiction.

The case $i-1 \in v_j \cap u_j$ is similar.
\end{proof}

\subsection{Duality between $\mMPR$ and $\MMPR$}\label{sec:dual}
Let $V$ be the space of infinite $\Z$-linear combinations of
elements $\{v_s \mid s \in S\}$ where $S$ is some indexing set.  As
in Section~\ref{sec:mMPR} call a linear functional $f: V \to \Z$
continuous if it respects infinite linear combinations, not just
finite ones. The set of all continuous linear functionals forms the
{\it continuous dual} $V^*_c$ of $V$ and is a free $\Z$-module with
basis $\{f_s \mid s \in S\}$ where $f_s$ is defined by $f_s(v_{s'})
= \delta_{ss'}$. The $\Z$-module $V$ is then the usual dual of
$V^*_c$.  Abusing notation slightly, we may simply say that $V$ and
$V^*_c$ are {\it continuous duals} with dual bases $\{v_s \mid s \in
S\}$ and $\{f_s \mid s \in S\}$ (even though $\{v_s\}$ may {\it not}
be a basis of $V$).  This notion of continuous duality makes sense
for bialgebras, and Hopf algebras, with distinguished bases.

There is a natural way to consider $\M$-permutations and
$\m$-permutations to be inverses of each other.  If $w = w_1\cdots
w_k \in\bS_n$ then $u = w^{-1}$ is the $\m$-permutation of $k$ with
length $n$ such that the $i$-s in $u$ occur in the positions
specified by $w_i$.  Clearly this gives rise to a bijection between
$\M$-permutations and $\m$-permutations.  For example, if $w$ is the
$\m$-permutation
$$[1,5,1,4,2,4,2,6,2,3]$$ then $w^{-1}$ is the $\M$-permutation
$$[(1,3),(5,7,9),(10),(4,6),2,8].$$  We note that the two
standardization operators are compatible in the following manner:
$\st(w) = v$ if and only if $\st(w^{-1}) = v^{-1}$, where we have
extended the inverse operation to words and set-compositions.
Furthermore, if $v \in S^\m_n$ then $v^{-1} \in \bS_{\ell(v)}$ and
$\ell(v^{-1}) = n$.

\begin{theorem}\label{thm:dual}
The bialgebras $\MMPR$ and $\mMPR$ are continuous duals of each
other with dual bases $\{w \mid w \in S^\m_\infty\}$ and $\{v \mid
v\in \bS_\infty\}$, where $S^\m_\infty$ is identified with
$\bS_\infty$ via $w \leftrightarrow w^{-1}$. The antipode
$S_{\MMPR}$ of $\MMPR$ induces an antipode $S_{\mMPR}$ of $\mMPR$
making $\MMPR$ and $\mMPR$ continuous duals as Hopf algebras.
\end{theorem}

\begin{proof}
The preceding comments on continuous duals allow us to prove the
theorem by simply comparing structure constants.  We first show that
product structure constants of $\mMPR$ are equal to coproduct
structure constants of $\MMPR$.  Let $u \in S^\m_n$, $v \in S^\m_m$
and $w \in S^\m_{m+n}$.  Since $*$ and $\vartriangle$ are
multiplicity free, we must show that $w$ occurs in $u * v$ if and
only if $u^{-1} \otimes v^{-1}$ occurs in $\vartriangle w^{-1}$.
Suppose that $w$ occurs in $u * v$.  If $w^{-1} = w'_1 w'_2 \cdots
w'_{n+m}$ then $u^{-1} = \st(w'_1 \cdots w'_n)$ and $v^{-1} =
\st(w'_{n+1} \cdots w'_{n+m})$.  The converse is also clear.
 Similarly, one checks that $w^{-1}$ occurs in $u^{-1}
\bullet v^{-1}$ if and only if $u \otimes v$ occurs in $\Delta w$.

Finally, we need to check that that the antipode $S_{\MMPR}:\MMPR
\to \MMPR$ induces a well defined map $S_{\mMPR} = S^*_{\MMPR}:\mMPR
\to \mMPR$ (here $S^*_{\MMPR}$ denotes the continuous transpose).
Let $w \in \bS_n$. Then by (\ref{eq:antipode}), $S_{\MMPR}(w)$ is a
linear combination of the basis elements $u \in \cup_{i \leq n}
\bS_i $. Thus the continuous transpose $S^*_{\MMPR}: \mMPR \to
\mMPR$ sends the basis element $x = u^{-1} \in S^\m_k$ to a
(possibly infinite) linear combination of the basis elements
$\{y\mid \ell(y) \geq \ell(x)\}$.  We need to show that this extends
via continuous linearity to a well-defined map for arbitrary
elements of $\mMPR$. But given a fixed $y \in S^\m_\infty$ there are
only finitely many $x \in S^\m_\infty$ with shorter length.  Thus
the coefficient of $y$ in $S^*_{\MMPR}(\sum_{x \in S^\m_\infty} a_x
\; x)$ is well-defined.
\end{proof}
%

\section{The Hopf algebra of Multi-noncommutative symmetric functions} \label{sec:MNSym}

Now we define a Hopf subalgebra of $\MMPR$ which will turn out to be
the dual Hopf algebra of $\mQSym$.  Let $\alpha$ be a composition of
$n$.  We say that $w \in \bS_n$ is of {\it type} $\type(w) = \alpha$
if $w^{-1} \in S^\m_\infty$ has descent set $\Des(w^{-1}) =
D(\alpha)$. Alternatively, the type of $w$ is $\C(D)$ where
$$
D = \{i \in [1,n-1] \mid i+1 \;\text{lies to the right of $i$ in
$w$}\}.
$$
Let $\R_{\alpha} \in \MMPR$ denote the sum of all $w \in \bS_n$ of
type $\alpha$.  For example, if $\alpha = (3,1)$ then
$$
\R_{\alpha} = [(1,4),2,3]+[1,(2,4),3]+[1,2,4,3]+[1,4,2,3]+[4,1,2,3].
$$
Let $\MNSym$ denote the subspace of $\MMPR$ spanned by the elements
$\R_\alpha$ as $\alpha$ varies over all compositions.  Clearly the
elements $\R_\alpha$ are independent.
\subsection{The product structure on $\MNSym$} \label{sec:prodmnsym}

For two compositions $\alpha = (\alpha_1,\ldots,\alpha_k)$ and
$\beta = (\beta_1,\ldots,\beta_l)$, recall the definitions of
$\alpha \vartriangleright \beta$ and $\alpha \vartriangleleft \beta$
from Section~\ref{sec:NSym}.  We now define $\alpha \cdot \beta =
(\alpha_1,\ldots,\alpha_k+\beta_1 -1,\ldots,\beta_l)$.  For example,
if $\alpha = (3,2,5,1)$ and $\beta = (4,2)$ then $\alpha
\vartriangleright \beta = (3,2,5,5,2)$, $\alpha \vartriangleleft
\beta = (3,2,5,1,4,2)$ and $\alpha \cdot \beta = (3,2,5,4,2)$.

\begin{prop} \label{prop:boxdotmnsym}
Let $\alpha$ be a composition of $m$ and $\beta$ be a composition of
$n$. We have
$$
\R_\alpha \bullet \R_\beta = \R_{\alpha \vartriangleright \beta} + \R_{\alpha\vartriangleleft \beta} +
\R_{\alpha \cdot \beta}.
$$
Thus $\MNSym$ is closed under the product $\bullet$ of $\MMPR$.
\end{prop}
\begin{proof}
Directly from the definition it is clear that $\R_\alpha \bullet
\R_\beta$ is a multiplicity free sum of certain $\M$-permutations $u
\in \bS_{m+n} \cup \bS_{m+n-1}$.  The type of $u \in \bS_{m+n}$ is
determined by $\type(u|_{[m]})$, $\type(\st(u|_{[m+1,m+n]}))$ and
whether $m$ lies in front of $m+1$.  Such a $u$ occurs in the
product $\R_\alpha \bullet \R_\beta$ if and only if $\type(u|_{[m]})
= \alpha$ and $\type(\st(u|_{[m+1,m+n]})) = \beta$.  The type of $u
\in \bS_{m+n-1}$ is determined by $\type(u|_{[m]})$ and
$\type(\st(u|_{[m,m+n-1]}))$ and occurs in $\R_\alpha \bullet
\R_\beta$ if and only if $\type(u|_{[m]}) = \alpha$ and
$\type(\st(u|_{[m,m+n-1]})) = \beta$. The three terms $\R_{\alpha
\vartriangleright \beta}$, $\R_{\alpha\vartriangleleft \beta}$, and
$\R_{\alpha \cdot \beta}$ correspond (in order) to the following
three possibilities for $u$: (a) $u\in \bS_{m+n}$ and $m$ occurs
before $m+1$, (b) $u\in \bS_{m+n}$ and $m$ occurs after $m+1$, and
(c) $u \in \bS_{m+n - 1}$.
\end{proof}

\begin{prop} \label{prop:MNSymrel}
The algebra $\MNSym$ is isomorphic to the free algebra on the
symbols $\R_\alpha$ for each composition $\alpha$, with relations
$\R_\alpha \bullet \R_\beta = \R_{\alpha \vartriangleright \beta} +
\R_{\alpha \vartriangleleft \beta} + \R_{\alpha \cdot \beta}$.
\end{prop}

\begin{proof}
We must show that the relation $\R_\alpha \bullet \R_\beta = \R_{\alpha
\vartriangleright \beta} + \R_{\alpha \vartriangleleft \beta} +
\R_{\alpha \cdot \beta}$ implies all possible relations amongst the
elements $\{R_{\alpha}\}$.  Assume we have a relation in $\MNSym$.
Using $\R_\alpha \bullet \R_\beta = \R_{\alpha \vartriangleright \beta}
+ \R_{\alpha \vartriangleleft \beta} + \R_{\alpha \cdot \beta}$ we may
make the relation linear.  If all terms cancel out we conclude that
the original relation is implied by our basic set of relations.
Otherwise we obtain a linear dependence amongst the $\R_\alpha$'s.
But this is impossible as the sets of $\M$-permutations involved in
each $\R_{\alpha}$ are disjoint, and the set of $\M$-permutations
forms a basis of $\MMPR$.
\end{proof}

Define $F_k = \R_{(k)} = [1,2,\ldots,k] \in \MNSym$.

\begin{prop} \label{prop:MNSymfree}
$\MNSym$ is freely generated over $\Z$ by $\{F_k \mid k \geq 1\}$ as
an algebra.
\end{prop}

\begin{proof}
First we show that any $\R_{\alpha}$ can be written as polynomial in
the $F_k$'s. The proof proceeds by induction on the number of parts
in $\alpha$. The base case $\alpha = (k)$ is clear.

Let $\alpha = (\alpha_1, \ldots, \alpha_l)$ with $l \geq 2$ and let
$\bar \alpha = (\alpha_1, \ldots, \alpha_{l-1})$. Using
Proposition~\ref{prop:boxdotmnsym} to write $ \R_{\alpha} = \R_{\bar
\alpha} F_{\alpha_l} - \R_{\alpha'} - \R_{\alpha''}$ where $\alpha'
= \bar \alpha \vartriangleright \alpha_l = (\alpha_1, \ldots,
\alpha_{l-2}, \alpha_{l-1} + \alpha_l)$ and $\alpha'' = \bar \alpha
\cdot \alpha_l = (\alpha_1, \ldots, \alpha_{l-2}, \alpha_{l-1} +
\alpha_l-1)$.  Since $\R_{\bar \alpha}$, $\R_{\alpha'}$ and
$\R_{\alpha''}$ all have less parts than $\alpha$ by the inductive
assumption we may assume that they can be expressed in terms of the
$F_k$'s; therefore $\R_\alpha$ can also be written in terms of the
$F_k$'s.

Now we want to show that $F_k$'s are algebraically independent.
Assume that is not so and we have a non-trivial polynomial relation
$r$. Pick the monomial $F_{k_1} \dotsc F_{k_m}$ in $r$ with largest
total size $k_1 + \dotsc + k_m$, and amongst those pick one with $m$
largest. Now expanding $r$ in terms of the $w_\alpha$ basis, we see
that $\R_{(k_1, \ldots, k_m)}$ cannot come from any other monomial in
$r$, giving a contradiction.
\end{proof}

\subsection{The duality between $\MNSym$ and $\mQSym$}

\begin{theorem}
The subspace $\MNSym$ is a Hopf subalgebra of $\MMPR$ (continuously)
dual to $\mQSym$, and the basis $\{\R_\alpha\} \subset \MNSym$ is
dual to $\{L_\alpha\} \subset \mQSym$.
\end{theorem}

\begin{proof}
From the definitions it is clear that the subspace $\MNSym \subset
\MMPR$ is the annihilator of the space $I$ of Lemma~\ref{lem:mQSym}.
Since $I$ is a biideal, it follows immediately from
Theorem~\ref{thm:dual} that $\MNSym$ is a closed under product and
coproduct.  Thus $\MNSym$ is a Hopf subalgebra dual to $\mQSym$. One
then checks that the pairing $\ip{.,.}:\MMPR \otimes \mMPR \to \Z$
satisfies $\ip{\R_\alpha,u+I} = \delta_{\alpha\beta}$ where $u$ is
any $\m$-permutation satisfying $\Des(u) = \beta$.
\end{proof}

\begin{prop}\label{prop:MNSymNSym}
The Hopf algebras $\MNSym$ and $\NSym$ are isomorphic via the map
$F_k \mapsto S_k$.
\end{prop}
\begin{proof}
By Proposition~\ref{prop:MNSymfree} the map $F_k \mapsto S_k$ is an
algebra isomorphism.  It remains to verify (see also
Section~\ref{sec:NSym}) that $\vartriangle F_k = \sum_{0 \leq j \leq
k} F_j \otimes F_{k - j}$, where $F_0 := 1$. This is immediate from
the definition $F_k = [1,2,\ldots, k]$.
\end{proof}

\section{The big Hopf algebra of Multi-symmetric functions} \label{sec:MSym}

We now define the big Hopf algebra of Multi-symmetric functions
$\MSym$.  As we will see $\MSym$ is isomorphic to $\Sym$ as a Hopf
algebra, but it is naturally equipped with a basis $\{g_{\lambda}\}$
distinct from the Schur basis $\{s_\lambda\}$.

\subsection{Reverse plane partitions}
Let $\lambda$ be a partition which we associate with its Young
diagram in English notation.  A {\it {reverse plane partition}} $T$
of shape $\lambda$ is a filling of the boxes in $\lambda$ with
positive integers so that the numbers are weakly increasing in rows
and columns.  We write $\sh(T) = \lambda$ for the shape of a reverse
plane partition $T$.  For a plane partition $T$, define $x^T :=
\prod_{i \in \P} x_i^{T(i)}$ where $T(i)$ is {\it {the number of
columns}} in $T$ containing one or more entries equal to $i$. Note
that this is not the usual weight assigned to a reverse plane
partition.  Now define the {\it dual stable Grothendieck
polynomials} $g_\lambda(x_1,x_2,\ldots) \in \Z[[x_1,x_2,\ldots]]$ by
$$g_{\lambda}(x_1,x_2,\ldots) = \sum_{\sh(T)=\lambda} x^T.$$
Similarly define the skew polynomials $g_{\lambda/\mu}$ by taking
the sum over reverse plane partitions of shape $\lambda/\mu$.

\begin{theorem} \label{thm:gcom}
The formal power series $g_{\lambda/\mu}(x_1,x_2,\ldots)$ are
symmetric functions.
\end{theorem}

\begin{example}
Using the definition one computes
$$g_{(2,1)} =
m_{(2,1)}+2m_{(1,1,1)}+m_{(2)}+m_{(1,1)} = s_{(2,1)}+s_{(2)}.$$ Here
the $m_{\lambda}$'s are monomial symmetric functions, see
\cite{EC2}.
\end{example}

We give two proofs of this fact. Note that the top homogenous
component of $g_{\lambda}$ is just the Schur function $s_{\lambda}$
and thus is symmetric.

\begin{remark}
The dual stable Grothendieck polynomials $g_{\lambda}$ are
implicitly studied by Lenart~\cite{Len}.  Shimozono and Zabrocki
give a determinantal formula for the $g_\lambda$ in \cite{SZ}.
\end{remark}

\subsection{Fomin-Greene operators again}
\label{sec:FGagain} Let $\Z\Lambda$ be the free $\Z$-module of
formal linear combinations of partitions.  For each $i \in \P$, we
define linear operators $u_i: \Z\Lambda \to \Z\Lambda$ called {\it
column adding operators} as follows:
$$u_i(\lambda) = \sum \mu,$$ where the sum is over all valid Young
diagrams $\mu$ obtained from $\lambda$ by adding several (at least
one) cells to the $i$-th column.  If no boxes can be added to the
$i$-th column of $\lambda$ then $u_i(\lambda) = 0$.

\begin{example}
If $\lambda = (4,3,3,1)$ then $u_6(\lambda) = 0$, $u_5(\lambda) =
(5,3,3,1)$, $u_4(\lambda) = (4,4,3,1)+(4,4,4,1)$, $u_3(\lambda) =
0$, $u_2(\lambda)=(4,3,3,2)$, $u_1(\lambda) =
(4,3,3,1,1)+(4,3,3,1,1,1)+\ldots$. Note that $u_1(\lambda)$ is
always an infinite expression.
\end{example}

In \cite{FG} Fomin and Greene prove the following statement.

\begin{lemma} {\cite[Lemma~3.1]{FG}} \label{lem:ecom}
Assume that a set of elements $\{u_i \mid i \in \Z\}$ of an
associative algebra satisfy the relations
\begin{align*}
u_i u_k u_j &= u_k u_i u_j & \mbox{for $i<j<k$}\\
u_j u_i u_k &= u_j u_k u_i & \mbox{for $i<j<k$} \\
u_j u_i (u_i+u_j)&=(u_i+u_j)u_j u_i & \mbox{for $i<j$.}
\end{align*}
Then the noncommutative analogs of elementary symmetric functions
$$
e_k(u_1, u_2, \ldots) = \sum_{a_1 > a_2 > \ldots > a_k}u_{a_1}u_{a_2}\dotsc u_{a_k}
$$
commute.
\end{lemma}
Note that the statement of Lemma~\ref{lem:ecom} is formal: the
$e_k(u_1,u_2,\ldots)$ should be considered as elements of an
appropriate completion.

\begin{lemma} \label{lem:relu}
The column adding operators $u_i: \Z\Lambda \to \Z\Lambda$ satisfy
the relations in Lemma~\ref{lem:ecom}.
\end{lemma}

\begin{proof}
The first two relations are straightforward since if $|i-k|>1$
operators $u_i$ and $u_k$ can easily be seen to commute.

Thus, it remains to argue that the third relation holds. Again, if
$|i-j|>1$ operators $u_i$ and $u_j$ commute and the relation
follows. Thus the only non-trivial case is $j = i+1$. Let $\mu$ be a
partition occurring with non-zero coefficient in $u_j u_i
(u_i+u_j)(\lambda)$.  Then $\mu$ can be obtained from $\lambda$ in
the following ways, where by $\lambda +_i a$ we denote operation of
adding $a \geq 1$ cells to $i$-th column of $\lambda$.

\begin{enumerate}
\item $\mu = \lambda +_j a' +_i b +_j a''$;
\item
$\mu = \lambda +_i b' +_i b'' +_j a$ and $\lambda +_i b' +_j a +_i
b''$ is a valid sequence;
\item
$\mu = \lambda +_i b' +_i b'' +_j a$ and $\lambda +_i b' +_j a +_i
b''$ is not a valid sequence.
\end{enumerate}
Thus $\mu$ differs from $\lambda$ in the $i$-th column by $b = b' +
b''$ squares and in the $(i+1)$-th column by $a = a' + a''$ squares.
We now show how the three cases above correspond to terms equal to
$\mu$ occurring in $(u_i+u_j)u_j u_i(\lambda)$.

In case (1), we rearrange the terms to get $\mu = \lambda +_i b +_j
a' +_j a''$, which must be a valid sequence.  In case (2), we biject
the expression $\lambda +_i b' +_i b'' +_j a$ with the valid
sequence $\lambda +_i b' +_j a +_i b''$.  Finally, in case (3), we
biject $\lambda +_i b' +_i b'' +_j a$ with the sequence $\lambda +_i
b +_j a' +_j a''$, where $a'<a$ is the maximal number such that
$\lambda +_i b' +_j a'$ is a valid sequence, and $a'' = a-a'$. Note
that $a'' \geq b' \geq 1$.

One can verify now that we get each possible summand of
$(u_i+u_j)u_j u_i(\lambda)$ equal to $\mu$ exactly once this way.
Indeed, we have the following three cases:

\begin{enumerate}
 \item $\mu = \lambda +_i b +_j a' +_j a''$ where $\mu = \lambda +_j a' +_i b +_j a''$ is a valid sequence;
 \item $\mu = \lambda +_i b' +_j a +_i b''$;
 \item $\mu = \lambda +_i b +_j a' +_j a''$ where $\mu = \lambda +_j a' +_i b +_j a''$ is not a valid sequence.
\end{enumerate}

These three cases correspond exactly to the three cases listed
before.

An illustration of the proof is given in Figure \ref{fig:gr1}.

\begin{figure}
\begin{center}
\input{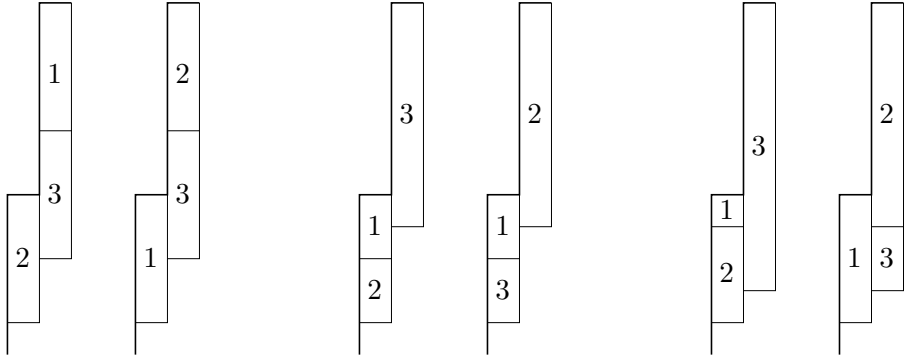}
\end{center}
\caption{The bijection between terms occurring in $u_j u_i
(u_i+u_j)(\lambda)$ and in $(u_i+u_j)u_j u_i(\lambda)$.  The three
possible cases are shown.}\label{fig:gr1}
\end{figure}
\end{proof}

\begin{proof}[Proof of Theorem~\ref{thm:gcom}]
Define the formal power series $A(x)$ with coefficients in operators
on $\Z\Lambda$ by $A(x) = \cdots (1+xu_2)(1+xu_1)$.
Lemma~\ref{lem:ecom} and Lemma~\ref{lem:relu} essentially says that
$A(x)A(y) = A(y)A(x)$ for two formal variables $x$ and $y$.  It is
immediate from the definition that
$$
g_{\lambda/\mu}(x_1,x_2,\ldots) = \ip{\cdots
A(x_2)A(x_1)\mu,\lambda}
$$
where $\ip{.,.}:\Z\Lambda \otimes \Z_{\rm fin}\Lambda \to \Z$ is the
pairing defined by $\ip{\nu,\rho} = \delta_{\rho\nu}$, and $\Z_{\rm
fin}\Lambda$ denotes the free $\Z$-module of {\it finite} linear
combinations of partitions.  Since the $A(x_i)$ commute,
$g_{\lambda/\mu}(x_1,x_2,\ldots)$ is a symmetric function in the
variables $x_i$.
\end{proof}

\begin{remark}  It would be interesting to compare the operators
$u_i$ of this section with the operators $v_i$ of
Section~\ref{sec:mSym}.  As we shall see the two kinds of operators
are dual in a sense which will be made clear in
Theorem~\ref{thm:MSymmSym}.
\end{remark}

\subsection{Schur decomposition of dual stable Grothendieck polynomials}
We give a direct bijection to establish a stronger version of
Theorem~\ref{thm:gcom}.  Namely, we describe an explicit rule for
the decomposition of the $g_{\lambda}$'s into Schur functions.

Given two partitions $\lambda$ and $\mu$ define the number
$f_{\lambda}^{\mu}$ as follows. If $\mu \nsubseteq \lambda$, we set
$f_{\lambda}^{\mu}=0$.  Otherwise, $f_{\lambda}^{\mu}$ is equal to
the number of {\it {elegant fillings}} of the skew shape
$\lambda/\mu$. A filling is elegant if it satisfies the following
two conditions:

\begin{enumerate}
\item
it is semistandard - that is, the numbers weakly increase in rows
and strictly in columns, and
\item the numbers in row $i$ lie in $[1, i-1]$.
\end{enumerate}
In particular there is no elegant filling of $\lambda/\mu$ if it
contains a cell in first row.  Elegant fillings were used previously
in~\cite{Len}.  An example of an elegant filling is given in
Figure~\ref{fig:gr2}.

\begin{figure}
\begin{center}
\input{gr2.pstex_t}
\end{center}
\caption{An elegant filling for $\lambda = (6,6,6,5,4,4)$ and $\mu =
(6,5,3,3,1)$.}\label{fig:gr2}
\end{figure}

\begin{theorem} \label{thm:gschur}
Let $\lambda$ be a partition.  Then $g_{\lambda}(x_1,x_2,\ldots) =
\sum_{\mu} f_{\lambda}^{\mu} s_{\mu}(x_1,x_2,\ldots)$.
\end{theorem}

\begin{proof}
We construct a weight preserving bijection between reverse plane
partitions $T$ of shape $\lambda$ and pairs $(S, U)$, where $S$ is a
semistandard tableau of shape $\mu$ and $U$ is an elegant filling of
shape $\lambda/\mu$.  Assume $\lambda$ has $m$ rows, and denote by
$T_i$ the $i$-th row of $T$.  More generally, let $T_{[i,j]}$ be the
reverse plane partition consisting of the part of $T$ between rows
$i$ and $j$ inclusively.  For each $i \in [1,m]$, define the {\it
{reduction}} $\tilde T_i$ of the row $T_i$ to be the sequence of
numbers obtained from $T_i$ by removing all entries equal to the
corresponding entries in $T_{i+1}$ immediately below.

To define the bijection we proceed recursively, defining a sequence
$(S_i,U_i)$ such that $S_i$ is a semistandard tableaux satisfying
$\sh(S_i) \subset (\lambda_i,\lambda_{i+1},\ldots,\lambda_m)$ and
$U_i$ is an elegant filling of shape
$(\lambda_i,\lambda_{i+1},\ldots,\lambda_m)/\sh(S_i)$.  For the
first step, set $S_m := T_m$ and $U_m$ to be the (empty) elegant
filling of $(\lambda_m)/(\lambda_m)$.  Assume now we have defined
$(S_{k+1}, U_{k+1})$.  We use the Robinson-Schensted-Knuth (RSK)
algorithm (see~\cite{EC2}) to insert the reduced row $\tilde T_k$
into $S_{k+1}$, obtaining $S_k$.  We verify by induction that
\begin{equation}\label{eq:firstrow}
\mbox{the first row of $S_i$ will always coincide with $T_i$.}
\end{equation}
Indeed, this is true for $S_m$, and the insertion of $\tilde T_i$
will always push out from the first row of $S_{i+1}$ elements of
$T_{i+1}$ which are strictly greater than elements of $T_i$
immediately above them inside $T$.

Now we describe how to obtain $U_k$.  First shift all the numbers in
$U_{k+1}$ one row down and simultaneously adding $1$ to each of
them, and consider the result as a partial filling of the cells of
the skew shape $\tau =
(\lambda_k,\lambda_{k+1},\ldots,\lambda_m)/\sh(S_{k+1})$.  Note that
the unfilled cells of $\tau$ form a horizontal strip $H$ of length
$\lambda_k$.  Now by well-known properties of RSK, the difference
$\sh(S_k)/\sh(S_{k+1})$ is a horizontal strip of size
$|\tilde{T}_k|\leq \lambda_k$, and this horizontal strip is
contained in $H$.  We obtain $U_k$ by placing a 1 in every cell that
lies in $H$ but not in $\sh(S_k)/\sh(S_{k+1})$, thus obtaining a
filling of shape
$(\lambda_k,\lambda_{k+1},\ldots,\lambda_m)/\sh(S_{k})$.  Since the
1's in $U_k$ form a horizontal strip, and by assumption $U_{k+1}$ is
semistandard, we conclude that $U_k$ is semistandard.  Also all the
entries in $U_{k+1}$ are moved one row down and incremented, and by
(\ref{eq:firstrow}) none of the new 1's are placed in the first row
so $U_k$ must be elegant, again assuming that $U_{k+1}$ was elegant.
Proceeding in this manner we obtain a pair $(S,U) :=(S_1, U_1)$.

Figure~\ref{fig:gr4} illustrates this direction of the bijection.

\begin{figure}
\begin{center}
\input{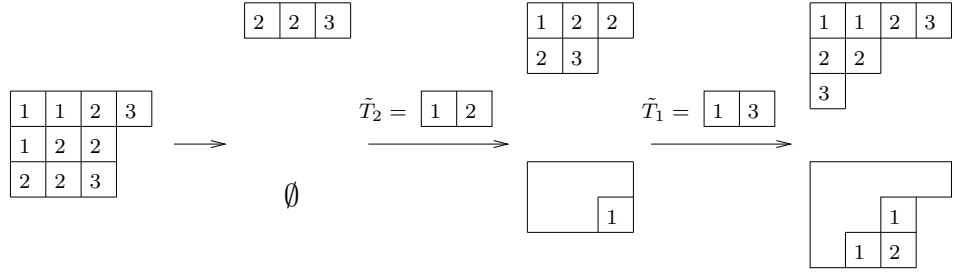}
\end{center}
\caption{The bijection in the proof of
Theorem~\ref{thm:gschur}.}\label{fig:gr4}
\end{figure}

To prove that the defined map is a bijection we describe the inverse
map. Start with a pair $(S, U) = (S_1, U_1)$ where $S$ is a
semistandard tableau of shape $\mu$ and $U$ is an elegant filling of
shape $\lambda/\mu$.  Construct recursively a sequence of pairs
$(S_i, U_i)$ such that $U_i$ has shape $\nu^{(i)}/\sh(S_i)$ for some
$\nu^{(i)}$, as follows. Assume $(S_k, U_k)$ has already been
constructed. Define the {\it {boundary}} of $S_k$ to be the set of
cells in $S_k$ directly below which there are no other cells. The
boundary is a horizontal strip with size equal to the size of the
first row of $S_k$. Define the {\it {active boundary}} of $S_k$ to
be the subset of cells of the boundary below which the cell of $U_k$
does not contain a 1. Note that active boundary is again always a
horizontal strip.  Now apply the inverse RSK algorithm to the active
boundary of $S_k$, producing a smaller semistandard tableau
$S_{k+1}$ and a non-decreasing sequence of numbers $\bar T_k$.  In
order to get $U_{k+1}$ remove the $1$'s from $U_k$, decrease all the
remaining numbers by $1$ and move them one row up. It is evident
that $U_{k+1}$ is an elegant filling if $U_k$ was.  By the choice of
active boundary it is also clear that the shape of $U_{k+1}$
``fits'' with the shape of $S_{k+1}$.

Now define $T$ by letting its $i$-th row $T_i$ equal the first row
of $S_i$.  By properties of the (inverse) RSK algorithm $T$ must be
a valid reverse plane partition.  Also observe that the sequences
$\bar T_i$ ejected during the construction are equal to the reduced
rows $\tilde T_i$ of $T$.  Indeed, by the nature of inverse RSK all
the elements of $\bar T_i$ were either bumped out from the first row
of $S_i$ by bigger numbers or belonged to the part of the active
boundary in the first row of $S_i$. This means that the numbers
$\bar T_i$ are exactly the elements of $T_i$ that are not equal to
the element of $T_{i+1}$ immediately below, which is the definition
of $\tilde T_i$.  This shows that defined map is indeed a bijection.

Since $S$ is obtained by inserting the reductions $\tilde T_i$ of
the rows of $T$, we have $x^T=x^S$.  Thus the bijection is
weight-preserving, completing the proof.
\end{proof}

\begin{example}
The decomposition $g_{(3,2,2)} =
s_{(3,2,2)}+2s_{(3,2,1)}+s_{(3,1,1)}+3s_{(3,2)}+2s_{(3,1)}+s_{(3)}$
corresponds to the following elegant fillings: \setcellsize{11}
$$
\tableau{{}&{}&{}\\{}&{}\\{}&{}} \ \tableau{{}&{}&{}\\{}&{}\\{}&{1}}
\ \tableau{{}&{}&{}\\{}&{}\\{}&{2}} \
\tableau{{}&{}&{}\\{}&{1}\\{}&{2}} \
\tableau{{}&{}&{}\\{}&{}\\{1}&{1}} \
\tableau{{}&{}&{}\\{}&{}\\{1}&{2}} \
\tableau{{}&{}&{}\\{}&{}\\{2}&{2}} \
\tableau{{}&{}&{}\\{}&{1}\\{1}&{2}} \
\tableau{{}&{}&{}\\{}&{1}\\{2}&{2}} \
\tableau{{}&{}&{}\\{1}&{1}\\{2}&{2}}
$$

\end{example}

Let $\MSym = \oplus_\lambda \Z g_{\lambda}$ be the free $\Z$-module
consisting of finite $\Z$-linear combinations of the $g_\lambda$.
\begin{prop}\label{prop:MSymSym}
The elements $g_{\lambda}$ form a basis for the ring of symmetric
functions.  Thus $\MSym \simeq \Sym$.
\end{prop}
\begin{proof}
By Theorem~\ref{thm:gschur}, the transition matrix between the basis
of Schur functions $\{s_\lambda\}$ and the set $\{g_{\lambda}\}$ is
upper triangular.
\end{proof}

$\MSym$ inherits from $\Sym$ a Hopf algebra structure.  While
isomorphic, they come with different distinguished bases:
$\{g_\lambda\}$ and $\{s_\lambda\}$.

\begin{prop}\label{prop:MSymfree}
$\MSym$ is freely generated by the set $\{g_{(n)} \mid n \geq 1\}$
as an algebra.
\end{prop}

\begin{proof}
Since $g_{(n)}=h_n$ are exactly the complete homogenous symmetric
functions, this follows from the well known fact that $\Sym =
\Z[h_1,\ldots,h_n]$.
\end{proof}

Let $\rho$ and $\tau$ be two skew shapes. Denote by $\rho
\vartriangleright \tau$ the skew shape obtained by attaching the two
so that lower leftmost cell of $\tau$ is directly to the right of
the upper rightmost cell of $\rho$.  Denote by $\rho \vartriangleleft \tau$ the
skew shape obtained by attaching the two shapes so that lower
leftmost cell of $\tau$ is directly above the upper rightmost cell
of $\rho$.  Finally, denote by $\rho \cdot \tau$ the skew
shape obtained by attaching the two shapes so that lower leftmost
cell of $\tau$ coincides with the upper rightmost cell of $\rho$.

Recall that in Section~\ref{sec:comsym} we have associated a ribbon
skew shape $r_\alpha$ to each composition $\alpha$.  The operations
$\alpha \vartriangleright \beta$, $\alpha \cdot \beta$ and $\alpha
\vartriangleleft \beta$ in Section~\ref{sec:prodmnsym} are
consistent with the ones introduced here:
$$
r_{\alpha \vartriangleright \beta} = r_{\alpha} \vartriangleright
r_{\beta} \ \ \ r_{\alpha \cdot \beta} = r_{\alpha} \cdot r_{\beta} \ \ \
r_{\alpha \vartriangleleft \beta} = r_{\alpha} \vartriangleleft r_{\beta}.
$$

\begin{lemma}\label{lem:gmult}
Let $\rho$ and $\tau$ be two skew shapes.  We have $g_{\rho}
g_{\tau} = g_{\rho \vartriangleright \tau} + g_{\rho \vartriangleleft \tau} - g_{\rho \cdot \tau}$.
\end{lemma}

\begin{figure}
\begin{center}
\input{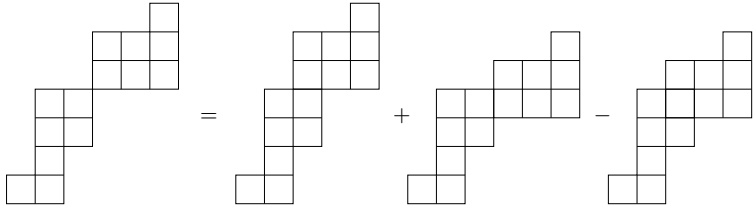}
\end{center}
\caption{The four shapes involved in
Lemma~\ref{lem:gmult}.}\label{fig:gr5}
\end{figure}

\begin{proof}
Let $R$ be a reverse plane partition of shape $\rho$ and $T$ be a
reverse plane partition of shape $\tau$. Let $R_{ur}$ be the label
of upper rightmost cell of $R$ and let $T_{ll}$ be the label of the
lower leftmost cell of $T$. If $R_{ur} \leq T_{ll}$ attach $R$ and
$T$ by putting $R_{ur}$ to the immediate left of $T_{ll}$ so that we
get a reverse plane partition $R \vartriangleright T$ of shape $\rho
\vartriangleright \tau$.  Note that $x^{R \vartriangleright T} = x^R
x^T$. If $R_{ur}
> T_{ll}$ attach the two so that we get a reverse plane partition $R \vartriangleleft T$ of
shape $\rho \vartriangleleft \tau$.  Note again that
$x^{R\vartriangleleft T} = x^R x^T$.

However, this map is not a bijection between pairs $(R,T)$ of
reverse plane partitions of shape $\rho$ or $\tau$ and a reverse
plane partition of shape $\rho \vartriangleright \tau$ or $\rho
\vartriangleleft \tau$.  There are additional reverse plane
partitions of shape $\rho \vartriangleleft \tau$ which cannot be
obtained in this way, namely, the ones where the box corresponding
to the lower leftmost box of $\tau$ has {\it {the same}} entry as
the box corresponding to upper rightmost cell of $\rho$.  Such
reverse plane partitions are in (weight-preserving) bijection with
reverse plane partitions of shape $\rho \cdot \tau$, which finishes
the proof.
\end{proof}

Figure~\ref{fig:gr5} illustrates Lemma~\ref{lem:gmult} by showing
the four shapes involved.  For each skew shape $\lambda/\mu$, define
the symmetric function $\tg_{\lambda/\mu}$ by
$$
\tg_{\lambda/\mu}(x_1,x_2,\ldots) =
(-1)^{|\lambda/\mu|}\,g_{\lambda/\mu}(-x_1,-x_2,\ldots).
$$
Thus $\tg_{\lambda/\mu}$ differs from $g_{\lambda/\mu}$ by a sign in
each homogeneous component.

\begin{theorem}
The map $F_k \mapsto g_{(k)} = \tg_{(k)}$ is a surjective Hopf
algebra morphism from $\MNSym$ to $\MSym$ which sends $\R_{\alpha}$
to $\tg_{r_\alpha}$ for each composition $\alpha$.  Thus, $\MSym$ is
the commutative image of $\MNSym$.
\end{theorem}
\begin{proof}
The first statement is immediate from Propositions
\ref{prop:MNSymfree}, \ref{prop:MNSymNSym}, \ref{prop:MSymSym} and
\ref{prop:MSymfree}.  Since $|\tau \cdot \rho| =
|\tau|+|\rho| - 1$, Lemma~\ref{lem:gmult} implies that $\tg_{\rho}
\tg_{\tau} = \tg_{\rho \vartriangleright \tau} + \tg_{\rho \cdot
\tau} + \tg_{\rho \vartriangleleft \tau}$.  For ribbon shapes, this
agrees with Proposition~\ref{prop:MNSymrel}, giving the statement of
the theorem.
\end{proof}

\subsection{Duality between $\MSym$ and $\mSym$}
In \cite{Len} Lenart proved the following theorem.
\begin{theorem} \cite[Theorem 2.8]{Len}
\label{thm:len} For a partition $\lambda$, one has $$ s_\lambda =
\sum_{\lambda \subset \mu} f^\lambda_\mu G_\mu
$$
where $f^\lambda_\mu$ is the number of elegant fillings of
$\mu/\lambda$.
\end{theorem}
Using Theorem \ref{thm:len}, we can relate the $g_\lambda$ to
$G_\lambda$.

\begin{theorem}\label{thm:MSymmSym}
The Hopf algebras $\MSym$ and $\mSym$ are continuously dual Hopf
algebras via the Hall inner product.  The bases $\{g_\lambda \mid
\lambda \in \Lambda\}$ and $\{G_\lambda \mid \lambda \in \Lambda\}$
are dual bases.  The bases $\{\tg_\lambda \mid \lambda \in
\Lambda\}$ and $\{\K_\lambda \mid \lambda \in \Lambda\}$ are dual
bases.
\end{theorem}
\begin{proof}
By Proposition~\ref{prop:MSymSym}, $\MSym$ is isomorphic to $\Sym$,
and by Proposition~\ref{prop:mSymSym}, $\mSym$ is isomorphic to the
completion of $\Sym$.  Since $\Sym$ is self dual under the Hall
inner product, $\MSym$ and $\mSym$ are continuously dual with this
pairing.

But Theorem~\ref{thm:len} with Theorem~\ref{thm:gschur}, one
immediately concludes that $\{G_\lambda \mid \lambda \in \Lambda\}$
and $\{g_\lambda \mid \lambda \in \Lambda\}$ are dual bases.  To
obtain the last statement we use the isomorphism $f(x_1,x_2,\ldots)
\mapsto f(-x_1,-x_2,\ldots)$ of $\Sym$.
\end{proof}

\subsection{$K$-homology of Grassmannians} \label{sec:khom}
Theorem~\ref{thm:MSymmSym} allows us to interpret the algebra
$\MSym$ as the $K$-homology of Grassmannians.  While $K$-homology
and $K$-cohomology are isomorphic for Grassmannians, we will find
that $\MSym$ is functorially covariant, like $K$-homology, while
$\mSym$ is contravariant in the corresponding sense.  We use the
notation introduced in Section~\ref{sec:mSym} and refer the reader
to \cite{B} for further details.

There is a pairing of $K^\circ \Gr(k,\c^n)$ and $K_\circ
\Gr(k,\c^n)$ obtained by the sequence $K^\circ \Gr(k,\c^n) \otimes
K_\circ \Gr(k,\c^n) \to K_\circ \Gr(k,\c^n) \to K_\circ(*) = \Z$,
where the first map is induced by taking tensor products and the
second map is the pushforward to a point.  If $\alpha \in
K^\circ\Gr(k,\c^n)$ and $\beta \in K_\circ \Gr(k,\c^n)$ we let
$(\alpha,\beta)$ denote this pairing.

Let $[\I_\lambda] \in K_\circ \Gr(k,\c^n)$ denote the class of the
{\it ideal sheaf} of the boundary of the Schubert variety
$X_\lambda$.  For $\lambda \subset (n-k)^k$, we let $\tilde \lambda
= (n - k - \lambda_k, \ldots, n- k - \lambda_1)$ denote the rotated
complement of $\lambda$ in the $(n-k)^k$ rectangle.  Buch shows in
\cite[p.30]{B} that the classes $[\I_\lambda]$ form a basis dual to
the classes $[\O_\lambda]$ of structure sheaves of Schubert
varieties.  More precisely, one has $([\O_\lambda],[\I_{\tilde
\mu}]) = \delta_{\lambda \mu}$.

Via Theorem~\ref{thm:Buch} we may identify the limit of this pairing
as $k,n \to \infty$ with the Hall inner product.  In this way one
may identify quotients of $\MSym$ with the $K$-homologies $K_\circ
\Gr(k,\c^n)$ of Grassmannians, as we now explain.  Again for
convenience we let $[\I_{\tilde \lambda}] = 0$ if $\lambda$ does not
fit in a $(n-k)^k$ rectangle.

\begin{theorem}
The map $\MSym \to K_\circ \Gr(k,\c^n)$ given by $g_\lambda \mapsto
[\I_{\tilde \lambda}]$ is a surjection.  It identifies the
comultiplication of $\MSym$ with the map $$\Delta_*: K_\circ
\Gr(k,\c^n) \to K_\circ \Gr(k,\c^n) \otimes K_\circ \Gr(k,\c^n)$$
induced by the diagonal embedding $\Delta:\Gr(k,\c^n) \to
\Gr(k,\c^n) \times \Gr(k,\c^n)$ and the multiplication of $\MSym$
with the maps
$$\phi_*:K_\circ \Gr(k_1,\c^{n_1}) \otimes K_\circ \Gr(k_2,\c^{n_2})
\to K_\circ \Gr(k_1+k_2,\c^{n_1+n_2})$$ induced by $\phi:
\Gr(k_1,\c^{n_1}) \times \Gr(k_2,\c^{n_2} \to
\Gr(k_1+k_2,\c^{n_1+n_2})$ (see discussion after
Theorem~\ref{thm:Buch}).
\end{theorem}
\begin{proof}
The first statement is clear from the definitions since the classes
$\{[\I_{\tilde \lambda}] \mid \lambda \subset (n-k)^k\}$ of the
ideal sheaves form a basis $K_\circ \Gr(k,\c^n)$.  We will check the
``comultiplication'' statement (the last statement is similar). Let
$X = \Gr(d,\c^n)$.  The product of two classes in $K^\circ X$ can be
calculated via the pullback $\Delta^*: K^\circ X \otimes K^\circ X
\to K^\circ X$ in $K$-theory: $\Delta^*([\O_\lambda] \otimes
[\O_\mu]) = [\O_\lambda].[\O_\mu]$.  By Theorem~\ref{thm:Buch}, the
coefficient of $G_\nu$ in $G_\lambda G_\mu$ is equal to the
coefficient of $[\O_\nu]$ in $[\O_\lambda].[\O_\mu]$.  This in turn
can be calculated via the projection formula as
$(\Delta^*([\O_\lambda] \otimes [\O_\mu]), [\I_{\tilde \nu}]) =
([\O_\lambda] \otimes [\O_\mu],\Delta_*[\I_{\tilde \nu}]) =
\sum_{\rho\tau} a^{\rho \tau}_{\nu}([\O_\lambda],[\I_{\tilde
\rho}])([\O_\mu],[\I_{\tilde \tau}]) = a^{\rho\tau}_\nu$ where
$\Delta_*[\I_{\tilde \nu}] = \sum_{\rho\tau}a^{\rho\tau}_\nu
[\I_{\tilde \rho}]\otimes[\I_{\tilde \tau}]$.  By
Theorem~\ref{thm:MSymmSym}, the product structure constants for
$\{G_\lambda \mid \lambda \in \Lambda\}$ agree with the coproduct
structure constants for $\{g_\lambda \mid \lambda \in \Lambda\}$. We
conclude that the comultiplication $\Delta_*$ of $K_\circ X$ agrees
with the comultiplication of $\MSym$.
\end{proof}

\subsection{Conjugate Fomin-Greene operators}
\label{sec:vst} Let $\{u_i \mid i \in \Z\}$ be a set of operators
satisfying Lemma~\ref{lem:ecom}. Recall that we have defined formal
power series
$$A(x) = \cdots(1+xu_1)(1+xu_0)(1+xu_{-1}) \cdots = \sum_{k \geq 0} e_k(u)x^k.$$
The $e_k(u)$'s commute and thus generate a
homomorphic image $\Sym(u)$ of the algebra of symmetric functions.
This allows to define $f(u)$ for any symmetric function $f$ as
follows: it is the image of $f \in \Sym$ under the map $\Sym \to
\Sym(u)$ given by $e_k \mapsto e_k(u)$.

Define the formal power series
$$B(x) = \cdots \frac{1}{1-xu_{-1}}\frac{1}{1-xu_0}\frac{1}{1-xu_1}
\cdots$$ where as usual $x$ is a formal variable commuting with all
the $u_i$.

\begin{lemma}
We have $B(x) = \sum_{k \geq 0} h_k(u) x^k$ where $h_0(u) = 1$.
\end{lemma}

\begin{proof}
For each $l \geq 0$ we have the well known identity $\sum_{k =
0}^{l} (-1)^k h_k e_{n-k}=0$.  From this one deduces that $A(x)
\tilde B(x) = 1$, where $\tilde B(x) = \sum_{k \geq 0} h_k(u) x^k$.
On the other hand $A(x)B(x) = 1$ also holds. This implies $\tilde
B(x) = B(x)$.
\end{proof}

\begin{lemma}
We have $$\cdots A(x_2) A(x_1) = \sum_\lambda s_\lambda(u)
s_{\lambda'}(x)$$ and $$\cdots B(x_2) B(x_1) = \sum_\lambda
s_\lambda(u) s_\lambda(x)$$ where the sums are over all partitions
$\lambda$.
\end{lemma}

\begin{proof}
Start with the usual Cauchy identity $$\prod_{i,j = 1}^{\infty}
(1+x_i y_j) = \sum_{\lambda} s_\lambda(y) s_{\lambda'}(x).$$ Group
the terms on the left hand side so that we get $$\prod_{j=1}^\infty
\sum_{k=0}^\infty e_k(y) x_j^k = \sum_{\lambda} s_\lambda(y)
s_{\lambda'}(x).$$ Now apply the transformation $\Sym \to \Sym(u)$
given by $e_k(y) \mapsto e_k(u)$ to both sides. We get exactly the
first equality. The proof of the second one is analogous.
\end{proof}

Now let us assume that the operators $u_i$ act on the space $\Z\Lambda$
of formal $\Z$-linear combinations of all partitions.  As before, we define
the inner product $\ip{\lambda,\mu} = \delta_{\lambda\mu}$.
Define $M_{\mu/\nu} = \ip{\cdots A(x_2) A(x_1) \cdot \nu,\mu}$ and
$N_{\mu/\nu} = \ip{\cdots B(x_2) B(x_1) \cdot \nu,\mu}$.  Since the
$u_i$ satisfy the conclusion of Lemma~\ref{lem:ecom}, both
$M_{\mu/\nu}$ and $N_{\mu/\nu}$ are symmetric functions in the
variables $x_1,x_2,\ldots$.

\begin{lemma} \label{lem:MN}
For a set of operators (or noncommutative variables) $\{u_i \mid i
\in \Z\}$ satisfying the conditions of Lemma~\ref{lem:ecom} we have
$\omega(M_{\mu/\nu}) = N_{\mu/\nu}$, where $\omega: \Sym \to \Sym$
is the algebra involution given by $e_k \mapsto h_k$.
\end{lemma}

\begin{proof}
We compute $$\omega(M_{\mu/\nu}) = \omega\left(\sum_\lambda
\ip{s_\lambda(u) \cdot \nu,\mu} s_{\lambda'}(x)\right) =
\sum_\lambda \ip{s_\lambda(u) \cdot \nu,\mu} s_\lambda(x) =
N_{\mu/\nu}.$$
\end{proof}

\subsection{Weak set-valued tableaux}
\label{sec:wst} Fix a partition $\nu$.  Now we specialize the
situation in Section~\ref{sec:vst} to the operators $\{v_i \mid i\in
\Z\}$ defined in Section \ref{sec:mSym} which act on
$\Z\Lambda_\nu$. The action of the product $1/(1-xv_1) 1/(1-xv_2)
\cdots$ on $\Z\Lambda_\nu$ then corresponds to the following {\it
{weak set-valued tableaux}}.

\begin{definition}
A {\it weak set-valued tableau} $T$ of shape $\lambda/\nu$ is a
filling of the boxes with finite non-empty multisets of positive
integers (thus, numbers in one box are not necessarily distinct) so
that
\begin{enumerate}
\item
the smallest number in each box is strictly bigger than the largest
number in the box directly to the left of it (if that box is present);
\item
the smallest number in each box is greater than or equal to the
largest number in the box directly above it (if that box
is present).
\end{enumerate}
For a {\it weak set-valued tableau} $T$, define $x^T$ to be
$\prod_{i \geq 1} x_i^{a_i}$ where $a_i$ is the number of
occurrences of the letter $i$ in $T$.
\end{definition}

This differs from the set-valued tableaux of Buch \cite{B} in two
ways: (a) the strict and weak inequalities have been swapped, and
(b) repeated numbers are allowed in each box.  The following weak
set-valued tableau $T$ has weight $x^T = x_1 x_2^3 x_3 x_4^3 x_5^2$.
\setcellsize{16}
$$
\tableau{{12}&{44} \\ {223}&{5} \\{45}}
$$


Let $J_{\lambda/\nu} = \sum_T x^T$ denote the weight generating
function of all weak set-valued tableaux $T$ of shape $\lambda/\nu$.
\begin{theorem}\label{thm:J}
We have $J_{\lambda/\nu}(x_1,x_2,\ldots) = \ip{\cdots B(x_2) B(x_1) \cdot \nu,\lambda}$ where
$$B(x) = \cdots \frac{1}{1-xv_{-1}}\frac{1}{1-xv_0}\frac{1}{1-xu_1}
\cdots.$$  In particular, $J_{\lambda/\nu}(x_1,x_2,\ldots)$ is a symmetric function in the variables $x_1,x_2,\ldots$.
\end{theorem}
\begin{proof}
The result is established in the same way as \cite[Theorem 3.1]{B}.
The multiple occurences of a single number in a box correspond to
the degree $2$ and higher terms of the expansion $1/(1-xv_i) = 1 +
xv_i + x^2 v_i^2 + \cdots$. The reversal of the order of operators
changed the strict and weak inequalities; or in other words, swapped
the notions of horizontal and vertical strips.
\end{proof}

The following is a direct consequence of (the proof of) Lemma
\ref{lem:MN}, Lemma~\ref{lem:K} and Theorem~\ref{thm:J}.

\begin{prop}
For any skew shape $\lambda/\nu$, we have $\omega(\K_{\lambda/\nu})
= J_{\lambda/\nu}$.
\end{prop}

\subsection{Valued-set tableaux}
\label{sec:vast}Now let $u_i:\Z\Lambda \to \Z\Lambda$ be the
operators defined in Section \ref{sec:FGagain}.

\begin{definition}
A {\it valued-set tableaux} $T$ of shape $\lambda/\mu$ is a filling
of the boxes of $\lambda/\mu$ with positive integers so that
\begin{enumerate}
\item the transpose of this filling of $T$ is a (usual) semistandard tableau, and
\item
we are provided with the additional information of a decomposition
of the shape into a disjoint union $\lambda/\mu = \sqcup A_j$  of
groups $A_j$ of boxes so that each $A_i$ is connected and completely
contained within a single column and all boxes in each $A_i$
contains the same number.
\end{enumerate}
For a {\it valued-set tableau} $T$, define $x^T$ to be
$\prod_{i \geq 1} x_i^{a_i}$ where $a_i$ is the number of
groups $A_j$ which contain the letter $i$.
\end{definition}

An example of a valued-set tableau can be seen in
Figure~\ref{fig:gr8}. The grouping of the boxes is shown by omitting
the edge separating the boxes.

\begin{figure}[h]
\begin{center}
\input{gr8.pstex_t}
\end{center}
\caption{An example of a valued-set tableau $T$ with shape
$(4,4,4,2,2,1)$.  We have $x^T=x_1^2 x_2^3 x_3^3 x_4 x_5^2
x_6$.}\label{fig:gr8}
\end{figure}

Let $j_{\lambda/\mu} = \sum_T x^T$ denote the generating function of
all valued-set tableaux of shape $\lambda/\mu$.

\begin{theorem}\label{thm:j}
We have $j_{\lambda/\nu}(x_1,x_2,\ldots) = \ip{\cdots B(x_2) B(x_1) \cdot \nu,\lambda}$ where
$$B(x) = \cdots \frac{1}{1-xu_{0}}\frac{1}{1-xu_1}\frac{1}{1-xu_2}
\cdots.$$  In particular, $j_{\lambda/\nu}(x_1,x_2,\ldots)$ is a
symmetric function in the variables $x_1,x_2,\ldots$.
\end{theorem}
\begin{proof}
The operator $B(x)$ acting on a partition $\nu$ adds a vertical
strip.  If $\mu/\nu$ is a vertical strip, the coefficient $x^k$ in
$\ip{B(x)\cdot \nu,\mu}$ is the number of ways to write each column
of $\mu/\nu$ as a disjoint union of non-empty groups of boxes, using
$k$ groups in total.  This recovers the definition of a valued-set
tableau.  The last statement follows from $B(x)B(y) = B(y)B(x)$.
\end{proof}

The following is a direct consequence of Lemma \ref{lem:MN} and Theorem~\ref{thm:J}.

\begin{prop}
We have $\omega(g_{\lambda/\mu}) = j_{\lambda/\mu}$.
\end{prop}

Note that since $\omega:\Sym \to \Sym$ is an algebra automorphism,
the  $K$-theory and $K$-homology of Grassmanians can also be
described in terms of $J_{\lambda}$'s and $j_{\lambda}$'s.

\end{document}